\documentstyle[amscd]{amsart}
\numberwithin{equation}{section}
\theoremstyle{plain}
\newtheorem{thm}{Theorem}[section]
\newtheorem{cor}[thm]{Corollary}
\newtheorem{lem}[thm]{Lemma}
\newtheorem{prop}[thm]{Proposition}

\begin{document}
\title{A certain synchronizing property of subshifts and flow equivalence}
\author{Kengo Matsumoto}
\address{ Department of Mathematics, 
Joetsu University of Education,
Joetsu 943-8512 Japan}
\email{kengo{@@}juen.ac.jp}
\begin{abstract}
We  will study  a certain synchronizing property of subshifts called $\lambda$-synchronization.
The $\lambda$-synchronizing subshifts form a large class of irreducible subshifts
containing irreducible sofic shifts.
We prove that the $\lambda$-synchronization is invariant under flow equivalence of subshifts.
The $\lambda$-synchronizing K-groups
and the $\lambda$-synchronizing Bowen-Franks groups
are studied and proved to be invariant under flow equivalence 
 of $\lambda$-synchronizing subshifts.
They are new flow equivalence invariants for $\lambda$-synchronizing subshifts.  
\end{abstract}
\maketitle
\def\Zp{{ {\Bbb Z}_+ }}
\def\U{{ {\cal U} }}
\def\S{{ {\cal S} }}
\def\M{{ {\cal M} }}
\def\P{{ {\cal P} }}
\def\Q{{ {\cal Q} }}
\def\G{{ {\cal G} }}

\def\LLL{{ {\frak L}^{\lambda(\Lambda)} }}
\def\LLTL{{ {\frak L}^{\lambda(\widetilde{\Lambda})} }}

\def\OFL{{ {\cal O}_{\frak L}}}
\def\OLL{{ {\cal O}_{\lambda(\Lambda)}  }}
\def\OLTL{{ {\cal O}_{\lambda(\widetilde{\Lambda})}  }}
\def\ALL{{ {\cal A}_{\lambda(\Lambda)}  }}
\def\ALTL{{ {\cal A}_{\lambda(\widetilde{\Lambda})}  }}
\def\DLL{{ {\cal D}_{\lambda(\Lambda)}  }}
\def\DLTL{{ {\cal D}_{\lambda(\widetilde{\Lambda})}  }}
\def\FL{{{\cal F}_{\frak L}}}
\def\FKL{{ {\cal F}_k^{l} }}
\def\Ext{{{\operatorname{Ext}}}}
\def\Im{{{\operatorname{Im}}}}
\def\Hom{{{\operatorname{Hom}}}}
\def\Ker{{{\operatorname{Ker}}}}
\def\Coker{{{\operatorname{Coker}}}}
\def\dim{{{\operatorname{dim}}}}
\def\id{{{\operatorname{id}}}}
\def\OLF{{{\cal O}_{{\frak L}^{Ch(D_F)}}}}
\def\OLN{{{\cal O}_{{\frak L}^{Ch(D_N)}}}}
\def\OLA{{{\cal O}_{{\frak L}^{Ch(D_A)}}}}
\def\LCHDA{{{{\frak L}^{Ch(D_A)}}}}
\def\LCHDF{{{{\frak L}^{Ch(D_F)}}}}
\def\LCHLA{{{{\frak L}^{Ch(\Lambda_A)}}}}
\def\LWA{{{{\frak L}^{W(\Lambda_A)}}}}


Keywords: $\lambda$-synchronizing subshifts,
flow equivalence, K-groups, Bowen-Franks groups.

Mathematics Subject Classification 2000:
Primary 37B10; Secondary 19K99, 46L80.

\section{Introduction}

 In \cite{KM2010},
 a certain synchronizing property called $\lambda$-synchronization has been introduced.
 The $\lambda$-synchronizing property  for a subshift $\Lambda$ is an equivalent property to
 the property D for the transpose of $\Lambda$,
which has been introduced by W. Krieger in \cite{Kr7}. 
As the  $\lambda$-synchronizing property  
is weaker than the  usual synchronizing property,
synchronizing subshifts are 
$\lambda$-synchronizing.
Hence irreducible sofic shifts are $\lambda$-synchronizing 
as well as Dyck shifts, $\beta$-shifts, Morse shifts, etc.
are $\lambda$-synchronizing.
Many irreducible subshifts have this property.
A $\lambda$-graph system is a labeled Bratteli diagram with
an additional structure called $\iota$-map (\cite{1999DocMath}). 
A finite directed labeled graph gives rise to a $\lambda$-graph system 
with stationary vertices
so that a sofic shift is presented by a $\lambda$-graph system 
with stationary vertices.
Not only sofic shifts but also all subshifts may be presented by 
$\lambda$-graph systems. 
There is a canonical method to 
construct a $\lambda$-graph system from an arbitrary subshift.
If a subshift is sofic,
the canonically constructed $\lambda$-graph system is one 
given by the left Krieger cover graph for the sofic shift. 
Hence the canonically constructed $\lambda$-graph system from a subshift
is regarded as a generalization of a left Krieger cover graph.
In \cite{KM2010},
a construction of $\lambda$-graph systems from
$\lambda$-synchronizing subshifts has been introduced.
If a $\lambda$-synchronizing subshift is sofic,
the constructed $\lambda$-graph system is one given by the left Fischer cover graph
for the sofic shift.
Hence the constructed $\lambda$-graph system
from a
$\lambda$-synchronizing subshift
is regarded as a generalization of a left Fischer cover graph.
It is called the canonical $\lambda$-synchronizing $\lambda$-graph system
for a $\lambda$-synchronizing subshift.

In this paper, we will first  
 characterize
the canonical $\lambda$-synchronizing $\lambda$-graph system in an intrinsic way,
and prove that it has a unique synchronizing property.
We will also prove that it is minimal 
in the sense that there exists no proper $\lambda$-graph subsystem that presents the subshift.
In \cite{KM2010}, it has been proved that 
the K-groups and the Bowen-Franks groups for the  $\lambda$-synchronizing $\lambda$-graph system are invariant under topological conjugacy,
so that they yield topological conjugacy invariants of $\lambda$-synchronizing subshifts.
In the second part of this paper,
we will prove that  $\lambda$-synchronization is invariant under flow equivalence of subshifts.
Furthermore we will prove  that the K-groups and 
the  Bowen-Franks groups for the canonical $\lambda$-synchronizing $\lambda$-graph systems
for $\lambda$-synchronizing subshifts are invariant 
under flow equivalence of $\lambda$-synchronizing subshifts.  
They are new nontrivial flow equivalence invariants for a large class of subshifts. 
In \cite{2001ETDS},
the author has extended the Bowen-Franks groups
to general subshifts.
The  Bowen-Franks groups are computed as the 
 Bowen-Franks groups for the nonnegative matrices of the left Krieger covers
  if the subshifts are sofic.
The Bowen-Franks groups for the $\lambda$-synchronizing  $\lambda$-graph system
are computed as the 
 Bowen-Franks groups for the nonnegative matrices of the left Fischer covers 
 if the subshifts are sofic.

Throughout the paper, 
we denote by 
$\Zp$ and ${\Bbb N}$ the set of nonnegative integers 
and the set of positive integers
respectively.

\medskip

{\it Acknowledgment:}
The author would like to thank Wolfgang Krieger
for his various discussions  and constant encouragements. 

\section{$\lambda$-synchronizing subshifts}
Let 
$\Sigma$
be a finite set
with its discrete topology. 
We call it an alphabet and each member of it a symbol.
Let $\Sigma^{\Bbb Z}$, $\Sigma^{\Bbb N}$ 
be the infinite product spaces 
$\prod_{i=-{\infty}}^{\infty}\Sigma_{i}$, 
$\prod_{i=1}^{\infty}\Sigma_i$ 
where 
$\Sigma_{i} = \Sigma$,
 endowed with the product topology respectively.
 The transformation $\sigma$ on $\Sigma^{\Bbb Z}$ 
given by 
$\sigma((x_i)_{i \in {\Bbb Z}}) = (x_{i+1})_{i\in {\Bbb Z}}$ 
is called the full shift.
 Let $\Lambda$ be a shift invariant closed subset of $\Sigma^{\Bbb Z}$ i.e. $\sigma(\Lambda) = \Lambda$.
  The topological dynamical system 
  $(\Lambda, \sigma\vert_{\Lambda})$
   is called a subshift and simply written as $\Lambda$.
We denote by
$X_{\Lambda} (\subset \Sigma^{\Bbb N})$
the set of all right one-sided sequences appearing in $\Lambda$.
We denote by $|\mu|$ 
the length $k$ for a word $\mu = \mu_1\cdots\mu_k, \, \mu_i \in \Sigma$.
For a natural number $l \in {\Bbb N}$,
we denote by $B_l(\Lambda)$ 
the set of all words appearing in some $(x_i)_{i \in {\Bbb Z}}$ of $\Lambda$ 
with length equal to
$l$.
Put
$B_*(\Lambda) = \sqcup_{l=0}^\infty B_l(\Lambda)$
where $B_0(\Lambda) =\{ \emptyset \}$ the empty word.
For a word
$\mu =\mu_1\cdots \mu_k \in B_*(\Lambda)$,
 a right infinite sequence
$x =(x_i)_{i \in {\Bbb N}} \in X_\Lambda$
and
$l \in \Zp$,
put
\begin{align*}
\Gamma_l^-(\mu) & = \{ \nu_1 \cdots \nu_l \in B_l(\Lambda) 
\mid \nu_1 \cdots \nu_l \mu_1 \cdots \mu_k \in B_*(\Lambda) \}, \\
\Gamma_l^-(x) & = \{ \nu_1\cdots \nu_l \in B_l(\Lambda) 
\mid (\nu_1,\dots, \nu_l, x_1,x_2, \dots) \in X_\Lambda \},\\
\Gamma_l^+(\mu) & = \{ \omega_1 \cdots \omega_l \in B_l(\Lambda) 
\mid \mu_1 \cdots \mu_k \omega_1 \cdots \omega_l \in B_*(\Lambda) \}, \\
\Gamma_\infty^+(\mu) & = \{ y \in X_\Lambda 
\mid \mu y \in X_\Lambda \}, \\
\intertext{ and}
\Gamma_*^-(\mu) & = \cup_{l=0}^\infty \Gamma_l^-(\mu), 
\qquad
\Gamma_*^+(\mu) = \cup_{l=0}^\infty \Gamma_l^+(\mu).
\end{align*}
A  word
$\mu =\mu_1\cdots \mu_k \in B_*(\Lambda)$
for 
$l \in \Zp$
is said to be $l$-{\it synchronizing}
if the equality
\begin{equation*}
\Gamma_l^-(\mu) = \Gamma_l^-(\mu \omega) 
\end{equation*}
holds for all 
$\omega \in \Gamma_*^+(\mu)$. 
Denote by
$S_l(\Lambda)$ the set of all $l$-synchronizing words of $\Lambda$.
It is easy to see that
a word 
 $\mu \in B_*(\Lambda)$
is $l$-synchrnonizing
if and only if 
$\Gamma_l^-(\mu) = \Gamma_l^-(\mu x)$
for all $x \in \Gamma_\infty^+(\mu)$.
Recall  
that an irreducible subshift $\Lambda$ is defined to be $\lambda$-{\it synchronizing}\
if 
for any $\eta \in B_l(\Lambda)$ and $k \ge l$
there exists $\nu \in S_k(\Lambda)$ 
such that 
$\eta \nu \in S_{k-l}(\Lambda)$ (see \cite{KM2010}).

Balnchard and Hansel have introduced a class of subshifts called synchronizing shifts
which contains the irreducible
sofic shifts (\cite{BH}).
Let $\Lambda$ be a subshift over $\Sigma$.
A word $\omega \in B_*(\Lambda)$
is said to be intrinsically synchronizing
if $\mu\omega, \omega \nu \in B_*(\Lambda)$ 
for $\mu, \nu \in B_*(\Lambda)$
implies 
$\mu \omega \nu \in B_*(\Lambda)$.
An irreducible subshift $\Lambda$ is said to be synchronizing 
if $\Lambda$ has an intrinsically synchronizing word.
\begin{prop}
A synchronizing  shift is $\lambda$-synchronizing.
Hence an irreducible sofic shift is $\lambda$-synchronizing.
\end{prop} 
\begin{pf}
Let $\omega$ be an intrinsically synchronizing word.
Then the words of the form
$\xi \omega \in B_*(\Lambda)$ for 
$\xi \in B_*(\Lambda)$
are intrinsically synchronizing by definition.
As for any intrinsically synchronizing word $\zeta$ of $\Lambda$ and 
$l \in \Zp$, we have
$$
\Gamma_l^-(\zeta) = \Gamma_l^-(\zeta \mu)
$$
for $\mu \in \Gamma_*^+(\zeta)$.
Hence
any intrinsically synchronizing word
is $l$-synchronizing for all $l \in {\Bbb N}$.
Since
$\Lambda$ is irreducible,
for $\eta \in B_l(\Lambda)$
and $k \ge l$,
there exists $\xi \in B_*(\Lambda)$ 
such that
$\eta \xi \omega \in B_*(\Lambda)$.
As
$\nu = \xi \omega$ is
intrinsically synchronizing and hence
$k$-synchronizing,
one has 
$\eta \nu \in S_{k-l}(\Lambda)$.
This means that $\Lambda$
is $\lambda$-synchronizing.
\end{pf}
There exists a concrete  example of an irreducible subshift that is not $\lambda$-synchronizing
(see \cite{KM2010}).

We note that it has been proved that property D is invariant under topological conjugacy 
(\cite{Kr7}). 
The following propostion is written  in \cite{KM2010}.
\begin{prop}[\cite{Kr7},\cite{KM2010}]
The
$\lambda$-synchronization is  invariant under topological
conjugacy of subshifts.
\end{prop}
\begin{pf}
Suppose that two subshifts 
$\Lambda$ over $\Sigma$ 
and 
$\Lambda'$ over $\Sigma'$ are bipartitely related
in the sense of \cite{Na}.
There exist alphabets $C, D$,
specifications
$\kappa: \Sigma \longrightarrow CD$
$\kappa': \Sigma' \longrightarrow DC$
 and
 a bipartite subshift $\widehat{\Lambda}$ 
 over $C \sqcup D$ 
 related to a bipartite conjugacy between 
 $\Lambda$ and $\Lambda'$.
 We may naturally  extend $\kappa$ and $ \kappa'$ 
to $B_*(\Lambda)$ and $B_*(\Lambda')$ 
respectively. 
Let us assume that $\Lambda$ is $\lambda$-synchronizing.
For $\eta' \in B_{l'}(\Lambda')$ and $k' \ge l'$
with
$\eta' =  \alpha'_1\cdots\alpha'_{l'}$
and
$\kappa'(\alpha'_i) = d_ic_i, i=1,\dots,l'$.
Take symbols
$c_0 \in C, d_{l'+1}\in D$ 
such that
$c_0 d_1 c_1 d_2 c_2 \cdots d_{l'}c_{l'}d_{l'+1} \in B_*(\widehat{\Lambda}).$
Set
$\eta = \kappa^{-1}(c_0 d_1 c_1 d_2 c_2 \cdots d_{l'}c_{l'}d_{l'+1})$ 
that belongs to
$B_{l'+1}(\Lambda')$.
Put
$l = l' +1$ and $k = k' +1$
and then we have
$\eta \in B_l(\Lambda)$ with
$k \ge l$.
Since $\Lambda$ is $\lambda$-synchronizing,
there exists $\nu \in S_k(\Lambda)$ such that 
$\eta \nu \in S_{k-l}(\Lambda)$.
Let
$\kappa(\nu) = 
\tilde{c}_1\tilde{d}_1\tilde{c}_2\tilde{d}_2 \cdots \tilde{c}_n\tilde{d}_n.
$
 Take
 $\tilde{c}_{n+1} \in C$ 
such that the word 
$$
c_0 d_1 c_1 d_2 c_2 \cdots d_{l'}c_{l'}d_{l'+1} 
\tilde{c}_1\tilde{d}_1\tilde{c}_2\tilde{d}_2 \cdots \tilde{c}_n\tilde{d}_n\tilde{c}_{n+1}
$$
is admissible in $\widehat{\Lambda}$.
Put
$$
\nu' = {\kappa'}^{-1}(d_{l'+1} 
\tilde{c}_1\tilde{d}_1\tilde{c}_2\tilde{d}_2 \cdots \tilde{c}_n\tilde{d}_n\tilde{c}_{n+1})
$$
that belongs to
$B_{n+1}(\Lambda')$.
As 
${\kappa'}(\nu') =  d_{l'+1} \kappa(\nu) \tilde{c}_{n+1}
$
 and
 $\nu \in S_k(\Lambda)$,
 one has
$\nu' \in S_{k-1}(\Lambda')$
so that 
$\nu' \in S_{k'}(\Lambda')$.
By the equality
$$
c_0 \kappa'(\eta') \kappa'(\nu') = \kappa(\eta)\kappa(\nu)\tilde{c}_{n+1},
$$
the word
$\eta' \nu'$ is admissible in $\Lambda'$,
which belongs to
$S_{k' - l'}(\Lambda')$.
Therefore $\Lambda'$ is $\lambda$-synchronizing.
\end{pf}

\section{$\lambda$-synchronizing $\lambda$-graph systems}

A $\lambda$-graph system 
is a graphical object presenting a subshift (\cite{1999DocMath}). 
It is a generalization of a finite labeled graph and has a close relation
to a construction of a certain class of $C^*$-algebras (\cite{2002DocMath}).  
Let ${\frak L} =(V,E,\lambda,\iota)$ be 
a $\lambda$-graph system 
 over $\Sigma$ with vertex set
$
V = \cup_{l \in \Zp} V_{l}
$
and  edge set
$
E = \cup_{l \in \Zp} E_{l,l+1}
$
with a labeling map
$\lambda: E \rightarrow \Sigma$, 
and that is supplied with  surjective maps
$
\iota( = \iota_{l,l+1}):V_{l+1} \rightarrow V_l
$
for
$
l \in  \Zp.
$
Here the vertex sets $V_{l},l \in \Zp$
are finite disjoint sets.   
Also  
$E_{l,l+1},l \in \Zp$
are finite disjoint sets.
An edge $e$ in $E_{l,l+1}$ has its source vertex $s(e)$ in $V_{l}$ 
and its terminal  vertex $t(e)$ 
in
$V_{l+1}$
respectively.
Every vertex in $V$ has a successor and  every 
vertex in $V_l$ for $l\in {\Bbb N}$ 
has a predecessor. 
It is then required that there exists an edge in $E_{l,l+1}$
with label $\alpha$ and its terminal is  $v \in V_{l+1}$
 if and only if 
 there exists an edge in $E_{l-1,l}$
with label $\alpha$ and its terminal is $\iota(v) \in V_{l}.$
For 
$u \in V_{l-1}$ and
$v \in V_{l+1},$
put
\begin{align*}
E^{\iota}_{l,l+1}(u, v)
& = \{e \in E_{l,l+1} \ | \ t(e) = v, \iota(s(e)) = u \},\\
E_{\iota}^{l-1,l}(u, v)
& = \{e \in E_{l-1,l} \ | \ s(e) = u, t(e) = \iota(v) \}.
\end{align*}
Then we require a bijective correspondence preserving their labels between 
$
E^{\iota}_{l,l+1}(u, v)
$
and
$
E_{\iota}^{l-1,l}(u, v)
$
for each pair of vertices
$u, v$.
We call this property  the local property of $\lambda$-graph system. 
We call an edge in $E$ a labeled edge and a finite sequence of connecting labeled edges a labeled path.
If a labeled path $\gamma$ labeled $\nu$
starts at a vertex $v$ in $V_l$
and ends at a vertex $u$ in $V_{l+n}$,
we say that $\nu$ leaves $v$
and write 
$s(\gamma)=v, t(\gamma) = u, \lambda(\gamma) = \nu.$  
We henceforth assume that ${\frak L}$ is left-resolving, 
which means that 
$t(e)\ne t(f)$ whenever $\lambda(e) = \lambda(f)$ for $e,f \in E$.
For a vertex
$v \in V_l$ denote by
$\Gamma_l^-(v)$ the predecessor set of $v$
which is defined by the set of  words  of length $l$
appearing as labeled paths from a vertex in $V_0$ to the vertex $v$.  
 ${\frak L}$ is said to be predecessor-separated if 
$\Gamma_l^-(v) \ne \Gamma_l^-(u)$
whenever $u, v\in V_l$  are distinct.
A subshift $\Lambda$ is said to be presented by 
a $\lambda$-graph system 
${\frak L}$
if the set of admissible words of $\Lambda$
coincides with the set of labeled paths appearing somewhere in ${\frak L}$.
Two $\lambda$-graph systems
${\frak L} =(V,E,\lambda,\iota)$ over $\Sigma$
and
${\frak L}' =(V',E',\lambda',\iota')$ over $\Sigma$
are said to be isomorphic if there exist
bijections
$\varPhi_V:V \longrightarrow V'$
and 
$\varPhi_E:E \longrightarrow E'$
satisfying
$\varPhi_V(V_l) = V'_l$
and 
$\varPhi_E(E_{l,l+1}) = E'_{l,l+1}$
such that
they give rise to a labeled graph isomorphism
compatible to $\iota$ and $ \iota'$.
We note that any essential finite directed labeled graph 
${\cal G} = ({\cal V}, {\cal E}, \lambda)$ over $\Sigma$
with vertex set ${\cal V}$, 
edge set ${\cal E}$ and labeling map $\lambda:{\cal E}\longrightarrow \Sigma$
 gives rise to a $\lambda$-graph system
${\frak L}_{\cal G} =(V,E,\lambda,\iota)$
by setting
$V_l ={\cal V}, E_{l,l+1} ={\cal E}, \iota = \id$
for all $l \in \Zp$
(cf.\cite{2002DocMath}).

Two points $x, y \in X_{\Lambda}$ are said to be 
$l$-{\it past equivalent},
written as
$x \sim_l y$,
if $\Gamma_l^-(x) = \Gamma_l^-(y)$.
For a fixed 
$l \in \Zp$, 
let 
$F_i^l, i= 1, 2,\dots, m(l)$ 
be the set  of all $l$-past equivalence classes of $X_{\Lambda}$
so that $X_{\Lambda}$ is a disjoint union of 
 $F_i^l, i= 1, 2,\dots, m(l)$. 
 Then the canonical $\lambda$-graph system
$
{\frak L}^\Lambda =(V^\Lambda,E^\Lambda,\lambda^\Lambda,\iota^\Lambda )
$
for $\Lambda$ is defined as follows
(\cite{1999DocMath}). 
The vertex set $V_l^\Lambda$ 
at level $l$
consist of the sets 
 $F_i^l,i=1,\dots,m(l)$.
 We write an edge with label $\alpha$ 
    from the vertex 
 $F_i^l \in V_l^\Lambda$ 
 to the vertex
 $F_{j}^{l+1} \in V_{l+1}^\Lambda$ 
 if
$ \alpha x \in F_i^l$
 for some 
 $x \in F_j^{l+1}$.
We denote by $E_{l,l+1}^\Lambda$ 
the set of all edges from $V_l^\Lambda$ 
to $V_{l+1}^\Lambda$.
There exists a natural map $\iota^{\Lambda}_{l,l+1}$ 
from $V_{l+1}^\Lambda$ 
to
$V_l^\Lambda$ 
by mapping $F_j^{l+1}$ to $F_i^l $
 when 
  $F_i^l $ contains  $F_j^{l+1}$.
Set
$V^\Lambda = \cup_{l \in \Zp} V_l^\Lambda $
 and
$E^\Lambda = \cup_{l \in \Zp} E_{l,l+1}^\Lambda$.
The labeling of edges is denoted by
$\lambda^\Lambda:E^\Lambda \rightarrow \Sigma$.
 The canonical $\lambda$-graph system
$
{\frak L}^\Lambda 
$
is left-resolving and predecessor-separated 
and it presents $\Lambda$.

A $\lambda$-graph system
\begin{equation*}
{\frak L}^{\lambda(\Lambda)} 
=(V^{\lambda(\Lambda)}, E^{\lambda(\Lambda)},
\lambda^{\lambda(\Lambda)},\iota^{\lambda(\Lambda)})
\end{equation*}
for a $\lambda$-synchronizing subshift $\Lambda$
has been introduced in \cite{KM2010}.
It is regarded as a left Fischer cover version for 
a $\lambda$-synchronizing subshift $\Lambda$
whereas 
the canonical $\lambda$-graph system  
$
{\frak L}^\Lambda 
$
is regarded as a left Krieger cover version for 
a subshift $\Lambda$.
For $\mu, \nu \in B_*(\Lambda)$,
if $\Gamma_l^-(\mu) = \Gamma_l^-(\nu)$,
we say that $\mu$ is $l$-past equivalent to $\nu$
and write it as  
$\mu\underset{l}{\sim}\nu$.
\begin{lem}[\cite{KM2010}]
Let $\Lambda$ be a $\lambda$-synchronizing subshift.
Then we have
\begin{enumerate}
\renewcommand{\labelenumi}{(\roman{enumi})}
\item 
For $l \in {\Bbb N}$ and $\eta \in B_*(\Lambda)$,
there exists
$\mu \in S_l(\Lambda)$
such that $\eta \in \Gamma_l^-(\mu)$. 
\item
For $\mu \in S_l(\Lambda)$,
there exists
$\mu' \in S_{l+1}(\Lambda)$
such that 
$\mu\underset{l}{\sim}\mu'$.
\item
For $\mu \in S_l(\Lambda)$,
there exist
$\beta \in \Sigma$
and
$\nu \in S_{l+1}(\Lambda)$
such that 
$\mu\underset{l}{\sim}\beta\nu$.
\end{enumerate}
\end{lem}
\begin{pf}
(i)
 The assertion is direct from definition of $\lambda$-synchronization.

(ii)
For $\mu \in S_l(\Lambda)$ with
$|\mu | = n$,
put 
$k = n + l+1$.
As $\Lambda$ is $\lambda$-synchronizing,
there exists $\nu \in S_k(\Lambda)$
such that
$\mu \nu \in S_{k - n}(\Lambda)$.
Put $\mu'=\mu \nu  \in S_{l+1}(\Lambda)$
so that 
$\mu\underset{l}{\sim}\mu'$.

(iii)
For $\mu \in S_l(\Lambda)$ with
$\mu = \mu_1 \cdots \mu_n$,
put
$ k = n +l$.
As $\Lambda$ is $\lambda$-synchronizing,
there exists $\omega \in S_k(\Lambda)$
such that
$\mu \omega \in S_{k - n}(\Lambda)$.
Set
$\beta = \mu_1$
and
$\nu = \mu_2 \cdots \mu_n \omega$.
Since $\omega \in S_k(\Lambda)$,
one has 
$\nu \in S_{k-(n-1)}(\Lambda)$
so that 
$\nu \in S_{l+1}(\Lambda)$
and 
$\mu\underset{l}{\sim}\beta \nu$.
\end{pf}
Let
$V_l^{\lambda(\Lambda)}$ be the $l$-past equivalence classes of
$S_l(\Lambda)$.
We denote by $[\mu]_l$ 
the equivalence class of $\mu \in S_l(\Lambda)$.
For
$\nu \in S_{l+1}(\Lambda)$ and $\alpha \in \Gamma_1^-(\nu)$,
define a labeled edge from
$[\alpha \nu]_l\in V_l^{\lambda(\Lambda)}$ 
to
$[\nu]_{l+1} \in V_{l+1}^{\lambda(\Lambda)}$ labeled $\alpha$.
The set of such labeled edges are denoted by $E^{\lambda(\Lambda)}_{l,l+1}$.
Since
$S_{l+1}(\Lambda) \subset S_l(\Lambda)$,
we have a natural map
$
[\mu]_{l+1} \in V_{l+1}^{\lambda(\Lambda)} 
\longrightarrow 
[\mu]_l \in V_l^{\lambda(\Lambda)} 
$
denoted by
$\iota^{\lambda(\Lambda)}_{l,l+1}$.
\begin{prop}[\cite{KM2010}]
${\frak L}^{\lambda(\Lambda)} 
=(V^{\lambda(\Lambda)}, E^{\lambda(\Lambda)},
\lambda^{\lambda(\Lambda)}, \iota^{\lambda(\Lambda)})
$
defines a $\lambda$-graph system that presents $\Lambda$.
\end{prop}
\begin{pf}
We will show that the local property of $\lambda$-graph system
holds.
For $[\mu]_l \in V_l^{\lambda(\Lambda)}$
and
$[\mu]_{l+2} \in V_{l+2}^{\lambda(\Lambda)}$
with
$\mu \in S_l(\Lambda),\nu \in S_{l+2}(\Lambda)$,
suppose that there exists a labeled edge 
from
$[\mu]_l$ to $[\nu]_{l+1}$
labeled $\alpha \in \Sigma$.
Hence
$\alpha \nu \underset{l}{\sim}\mu$.
There exist an edge from
$[\alpha \nu]_{l+1}$ 
to
$[\nu]_{l+2}$
labeled $\alpha$
and an $\iota$-map from
$[\alpha \nu]_{l+1}$ to 
$[\alpha \nu]_l$.
On the other hand,
suppose that there exist an
$\iota$-map from
$[\omega]_{l+1}$
to
$[\mu]_l$
and an edge from
$[\omega]_{l+1}$
to
$[\nu]_{l+2}$
labeled $\alpha$
so that
$\omega\underset{l+1}{\sim}\alpha\nu$.
Since
$\iota^{\lambda(\Lambda)}([\alpha\nu]_{l+1})= 
[\alpha\nu]_l$,
one has 
$
\mu\underset{l}{\sim}
\alpha\nu$.
Hence there exists an edge from
$[\mu]_l$
to
$[\nu]_{l+1}$
labeled $\alpha$.
Therefore the local property of $\lambda$-graph system holds.
By definition of $\lambda$-synchronization of $\Lambda$,
an admissible word in $\Lambda$ appears in $\LLL$
as a labeled edge.
Hence $\LLL$ presents $\Lambda$.
\end{pf}
We call
$
{\frak L}^{\lambda(\Lambda)} 
$
the canonical
$\lambda$-synchronizing $\lambda$-graph system of $\Lambda$.
It is direct to see that 
$
{\frak L}^{\lambda(\Lambda)} 
$
is left-resolving and predecessor-separated.

Let ${\frak L} =(V,E,\lambda,\iota)$ 
be a $\lambda$-graph system over $\Sigma$
that presents a subshift $\Lambda$.
\noindent
{\bf Definition.} 
A $\lambda$-graph system
${\frak L}' =(V', E', \lambda', \iota')$
over $\Sigma'$
is called a  $\lambda$-graph subsystem of ${\frak L}$
if 
$\Sigma^\prime \subset \Sigma$
and
the following conditions hold for $l \in \Zp$:
\begin{equation*}
 V^{\prime}_l \subset V_l, 
\quad
E^{\prime}_{l,l+1} \subset E_{l,l+1},
\quad
\lambda^{\prime}_{l,l+1} = \lambda_{l,l+1} |_{E^{\prime}_{l,l+1}},
\quad
\iota^{\prime}_{l,l+1} = \iota_{l,l+1} |_{E^{\prime}_{l,l+1}}.
\end{equation*}
\begin{cor}
$
{\frak L}^{\lambda(\Lambda)} 
$
is a $\lambda$-graph subsystem
of 
$
{\frak L}^{\Lambda}. 
$
\end{cor}
\begin{pf}
Let
${\frak L}^{\Lambda} 
=(V^\Lambda, E^\Lambda,\lambda^\Lambda,\iota^\Lambda)
$
be the canonical $\lambda$-graph system
for
$\Lambda$.
An $l$-synchronizing word
$\mu \in S_l(\Lambda)$
satisfies
$
\Gamma_l^-(\mu) 
= 
\Gamma_l^-(\mu x)
$ 
for
$x \in \Gamma_\infty^+(\mu)$.
Hence
$[\mu]_l$ naturally defines  a vertex of $V_l^\Lambda$
so that
the vertex set
$V^{\lambda(\Lambda)}_l$ is regarded as a subset of $V_l^\Lambda$.
Similarly
the edge set
$E_{l,l+1}^{\lambda(\Lambda)}$ 
is regarded as a subset of $E_{l,l+1}^\Lambda$.
The $\iota$ map $\iota^{\lambda(\Lambda)}$ 
of ${\frak L}^{\lambda(\Lambda)}$
is obtained by restricting 
the $\iota$-map $\iota^{\Lambda}$
of ${\frak L}^\Lambda$.
Therefore    
$
{\frak L}^{\lambda(\Lambda)} 
$
is a $\lambda$-graph subsystem
of ${\frak L}^\Lambda$.
\end{pf}

We will characterize 
the canonical $\lambda$-synchronizing $\lambda$-graph system in an intrinsic way.
Let ${\frak L} = (V, E,\lambda,\iota)$
be a left-resolving, predecessor-separated $\lambda$-graph system 
over $\Sigma$ that presents a subshift $\Lambda$.
Denote by $\{ v_1^l,\dots,v_{m(l)}^l \}$
the vertex set $V_l$ at level $l$.
For an admissible word $\nu \in B_n(\Lambda)$
and a vertex $v_i^l \in V_l$,
we say that $v_i^l$ {\it launches}\, $\nu$
if the following two conditions hold:
\begin{enumerate}
\renewcommand{\labelenumi}{(\roman{enumi})}
\item
There exists a path labeled $\nu$ in ${\frak L}$ 
leaving  $v_i^l$
and ending at a vertex in $V_{l+n}$.
\item 
The word $\nu$
does not leave any other vertex in $V_l$
than $v_i^l$.
\end{enumerate}
The vertex $v_i^l$ is called the {\it launching vertex}\, for $\nu$.

\noindent
{\bf Definition.}
A $\lambda$-graph system  ${\frak L}$
is said to be $\lambda$-{\it synchronizing}\,
if for any $l \in \Zp$
and any vertex $v_i^l \in V_l$,
there exists a word $\nu \in B_*(\Lambda)$
such that 
$v_i^l$ launches $\nu$.
We set
$$
S_{v_i^l}(\Lambda)
=\{ \nu \in B_*(\Lambda) \mid v_i^l \text{ launches } \nu \}.
$$
\begin{lem}
Keep the above notations.
Assume that ${\frak L}= (V, E,\lambda,\iota)$ 
is $\lambda$-synchronizing.
Then we have
\begin{enumerate}
\renewcommand{\labelenumi}{(\roman{enumi})}
\item
$\sqcup_{i=1}^{m(l)} S_{v_i^l}(\Lambda) = S_l(\Lambda)$.
\item
The $l$-past equivalence classes of $S_l(\Lambda)$
is $S_{v_i^l}(\Lambda), i=1,\dots,m(l)$.
\item
For any $l$-synchronizing word $w \in S_l(\Lambda)$,
there exists a vertex $v_{i(\omega)}^l \in V_l$
such that $v_{i(\omega)}^l$ launches $\omega$ 
and
$
\Gamma_l^-(\omega)= \Gamma_l^-(v_{i(\omega)}^l).
$
\end{enumerate}
\end{lem}
\begin{pf}
(i)
The inclusion relation
$ S_{v_i^l}(\Lambda) \subset S_l(\Lambda)$
is obvious.
We will show that
$\sqcup_{i=1}^{m(l)} S_{v_i^l}(\Lambda) \supset S_l(\Lambda)$.
For
$\mu \in S_l(\Lambda)$,
suppose that $\mu$ leaves  two vertices
$v_i^l,v_j^l\in V_l,$
for 
$i\ne j$.
Since ${\frak L}$ is left-resolving,
there exist two distinct terminal vertices 
$v_{i'}^{l+n}, v_{j'}^{l+n}\in V_{l+n}$
with
$i' \ne j'$
of the labeled paths labeled $\mu$.
As ${\frak L}$ is $\lambda$-synchronizing,
for the vertices  
$v_{i'}^{l+n}, v_{j'}^{l+n}\in V_{l+n}$
there exist
admissible words
$\nu(i'), \nu(j') \in B_*(\Lambda)$
such that 
$v_{i'}^{l+n}$ launches $\nu(i')$
and
$v_{j'}^{l+n}$ launches $\nu(j')$.
As
$\nu(i')$ does not leave any other vertex than $v_{i'}^{l+n}$
in $V_{l+n}$,
one has
$\nu(i')\ne \nu(j')$.
One knows then that
\begin{equation*}
\Gamma_l^-(v_i^l) = \Gamma_l^-(\mu \nu(i')),
\qquad
\Gamma_l^-(v_j^l) = \Gamma_l^-(\mu \nu(j')).
\end{equation*}
Since
${\frak L}$ is predecessor-separated,
one has
\begin{equation*}
 \Gamma_l^-(\mu \nu(i'))
 \ne
 \Gamma_l^-(\mu \nu(j'))
\end{equation*}
and a contradiction to the hypothesis that
$\mu$ is a $l$-synchronizing word.

(ii)
For $\mu \in S_{v_i^l}(\Lambda)$,
$\nu \in S_{v_j^l}(\Lambda)$
with $i\ne j$,
one has 
\begin{equation*}
\Gamma_l^-(\mu) = \Gamma_l^-(v_i^l) 
\ne 
\Gamma_l^-(v_j^l) = \Gamma_l^-(\nu)
\end{equation*}
and hence
$[\mu]_l \ne [\nu]_l$.
Conversely for 
$\omega,\zeta \in S_l(\Lambda)$
with
$[\omega]_l \ne [\zeta]_l$,
as in the above discussion,
there uniquely exists a vertex 
$v_{i(\omega)}^l \in V_l$
such that 
$\omega$ leaves 
 $v_{i(\omega)}^l$.
This means that
$v_{i(\omega)}^l$ launches $\omega$.
Similarly,
there uniquely exists a vertex 
$v_{i(\zeta)}^l \in V_l$
such that 
$v_{i(\zeta)}^l$ launches $\zeta$.
Since
\begin{equation*}
 \Gamma_l^-(v_{i(\omega)}^l)= \Gamma_l^-(\omega),
 \qquad 
\Gamma_l^-(v_{i(\zeta)}^l) = \Gamma_l^-(\zeta)
\end{equation*}
and 
$\Gamma_l^-(\omega)
\ne
 \Gamma_l^-(\zeta),
$
we have
$i(\omega) \ne i(\zeta)$.

(iii)
The assertion is now clear from the above discussions.
\end{pf}

\noindent
{\bf Definition.}
A $\lambda$-graph system
${\frak L} = (V,E,\lambda,\iota)$ 
is said to be
$\iota$-{\it irreducible}\,
if for any two vertices $v,u \in V_l$ and 
a labeled path 
$\gamma$ starting at $u$,
there exist a labeled path 
from $v$ to a vertex $u'\in V_{l+n}$
such that
$\iota^n(u') = u$ ,
and a labeled path 
$\gamma'$ starting at $u'$
such that
$\iota^n(t(\gamma')) = t(\gamma)$
and
$\lambda(\gamma') = \lambda(\gamma)$,
where 
$t(\gamma'), t(\gamma)$
denote the terminal vertices of
$\gamma', \gamma$
respectively and
$\lambda(\gamma'),\lambda(\gamma)$
the words 
labeled  by 
$\gamma', \gamma$
respectively. 
A finite directed labeled graph ${\cal G}$
is irredsucible as a directed graph if and only if 
the $\lambda$-graph system
${\frak L}_{\cal G}$ is $\iota$-irreducible.

\begin{lem}
Let
${\frak L} = (V,E,\lambda,\iota)$
be a $\lambda$-graph system
that presents a subshift $\Lambda$.
If ${\frak L}$ is $\iota$-irreducible,
then
$\Lambda$ is irreducible.
\end{lem}
\begin{pf}
For $\mu,\nu \in B_*(\Lambda)$,
put
$k = |\mu|, l=|\nu|$.
Take a labeled  path labeled $\nu$
from a vertex in $V_0$ to a vertex $v$ in $V_l$.
Take a labeled path labeled $\mu$
from a vertex $u$ in $V_l$ to a vertex in $V_{l+k}$.
Since ${\frak L}$ is $\iota$-irreducible,
there exists a labeled path $\pi$ from
$v$ to a vertex $u'$ 
such that the word $\mu$ leaves $u'$.
Denote by
$\omega$ the word of the path $\pi$.
 We then have a labeled path which presents 
 the word
 $\nu \omega \mu \in B_*(\Lambda)$.
 Hence $\Lambda$ is irreducible.
\end{pf}
Conversely,
we have
\begin{lem}
Assume that
${\frak L} = (V,E,\lambda,\iota)$
is $\lambda$-synchronizing.
If 
$\Lambda$ is irreducible,
then
${\frak L}$ is $\iota$-irreducible.
\end{lem}
\begin{pf}
For two vertices $v,u \in V_l$ and 
a labeled path 
$\gamma$ starting at $u$,
put
$\mu = \lambda(\gamma)$
and 
$k = |\mu|$.
Let 
$w \in V_{l+k}$ 
denote the terminal vertex $t(\gamma)$
of $\gamma$.
Since ${\frak L}$ is $\lambda$-synchronizing,
$w$ is a launching vertex for a word $\eta \in B_*(\Lambda)$.
Similarly
$v$ is a launching vertex for a word $\zeta \in B_*(\Lambda)$.
Put
$n = |\zeta|$.
By the hypothesis that
$\Lambda$ is irreducible,
there exists a word $\xi\in B_*(\Lambda)$
such that 
$\zeta \xi \mu \eta \in B_*(\Lambda)$.
Since the word $\zeta$ must leave $v$ in $V_l$,
any labeled path labeled $\zeta \xi \mu \eta$
must leave $v$ in $V_l$.
Put
$m = |\xi|$.
Let $u' \in V_{l+n+m}$ be a vertex in $V_{l+n+m}$
at which the word $\zeta \xi$ ends and
the word $\mu \eta$ starts.
Denote by $\gamma'$ the labeled path labeled $\mu$ starting at $u'$.
Let $w'$ be the terminal vertex of $\gamma'$.
By the local property of $\lambda$-graph system,
there exists a labeled path labeled $\mu$
which starts at
$\iota^{n+m}(u') \in V_l$ 
and ends at 
$\iota^{n+m}(w') \in V_{l+k}$,
and there exists a labeled path
labeled $\eta$ starting at 
$\iota^{n+m}(w')$.
As $w$ is a launching vertex for $\eta$
and ${\frak L}$ is left-resolving,
we have
$\iota^{n+m}(w') = w$ 
and
$\iota^{n+m}(u') = u$.
This implies that 
${\frak L}$ is $\iota$-irreducible.
\end{pf}
Therefore we have
\begin{prop}
Let
${\frak L}$ be a $\lambda$-synchronizing $\lambda$-graph system that presents a subshift $\Lambda$.
Then $\Lambda$ is irreducible if and only if ${\frak L}$ is $\iota$-irreducible.
\end{prop}

\begin{prop}
A subshift $\Lambda$ 
is $\lambda$-synchronizing
if and only if
ther exists a left-resolving, predecessor-separated,
$\iota$-irreducible,
$\lambda$-synchronizing $\lambda$-graph system that presents $\Lambda$.
\end{prop}
\begin{pf}
Let ${\frak L}$ be 
a left-resolving, predecessor-separated, $\iota$-irreducible,
$\lambda$-synchronizing $\lambda$-graph system that presents $\Lambda$.
For 
any $\eta \in B_k(\Lambda)$ and $l \in {\Bbb N}$
with
$k \le l$,
there exists a terminal vertex $v_i^l \in V_l$ 
of a path labeled  $\eta$.
$\lambda$-synchronization of ${\frak L}$ 
implies  
$S_{v_i^l}(\Lambda) \ne \emptyset$.
Take  
a word
$\mu \in S_{v_i^l}(\Lambda)$.
As 
$S_{v_i^l}(\Lambda) \subset S_l(\Lambda)$,
we have
$\mu \in S_l(\Lambda)$
and 
$\eta \mu \in B_*(\Lambda)$.
By the previous lemma,
$\Lambda$ is irreducible.
Hence  
$\Lambda$ is $\lambda$-synchronizing.

Conversely suppose that 
$\Lambda$ is $\lambda$-synchronizing.
The canonical $\lambda$-synchronizing $\lambda$-graph system
${\frak L}^{\lambda(\Lambda)}$ for $\Lambda$
is a left-resolving, predecessor-separated,
$\lambda$-synchronizing $\lambda$-graph system that presents $\Lambda$.
As $\Lambda$ is irreducible,
the $\lambda$-synchronization of ${\frak L}$ implies
that 
${\frak L}^{\lambda(\Lambda)}$
is $\iota$-irreducible.
\end{pf}

\begin{thm}
For a $\lambda$-synchronizing
subshift $\Lambda$,
there uniquely exists a left-resolving, predecessor-separated,
$\iota$-irreducible,
$\lambda$-synchronizing $\lambda$-graph system that presents $\Lambda$.
The unique 
$\lambda$-synchronizing $\lambda$-graph system
is the canonical $\lambda$-synchronizing $\lambda$-graph system
${\frak L}^{\lambda(\Lambda)}$ for $\Lambda$.
\end{thm}
\begin{pf}
We will prove that 
the canonical $\lambda$-synchronizing $\lambda$-graph system
${\frak L}^{\lambda(\Lambda)}$ for $\Lambda$
is a unique 
 left-resolving, predecessor-separated, $\iota$-irreducible,
$\lambda$-synchronizing $\lambda$-graph system that presents $\Lambda$.
Let ${\frak L} =(V, E, \lambda, \iota)$
be a $\lambda$-graph system satisfying these properties.
For $u \in V_l^{\lambda(\Lambda)}$,
take a word
$\mu(u) \in S_l(\Lambda)$
such that 
$u = [\mu(u)]_l \in V^{\lambda(\Lambda)}_l$.
Since
${\frak L}$ is $\lambda$-synchronizing,
there  exists a unique vertex
$v_{\mu(u)}^l \in V_l$ that launches $\mu(u)$ and satisfies
$\Gamma_l^-(\mu(u)) = \Gamma_l^-(v_{\mu(u)}^l)$.
Define
$\varPhi_V: V_l^{\lambda(\Lambda)}\longrightarrow V_l$
by 
$\varPhi_V(u) =  v_{\mu(u)}^l$.
We will show that 
$\varPhi_V$ 
yields 
an isomorphism between ${\frak L}$ and ${\frak L}^{\lambda(\Lambda)}$
as in the following way.

1. Well-definedness of $\varPhi_V$ :
Let $\mu'(u) \in S_l(\Lambda)$ be another word such as 
$u =[\mu'(u)]_l$.
One then has 
$\Gamma_l^-(\mu(u)) = \Gamma_l^-(\mu'(u))$.
Since ${\frak L}$ is predecessor-separated,
one sees 
$v_{\mu(u)}^l = v_{\mu'(u)}^l$ in $V_l$.

2. $\iota\circ \varPhi_V = \varPhi_V\circ \iota^{\lambda(\Lambda)}$ :
For $w \in V_{l+1}^{\lambda(\Lambda)}$,
put
$w' = \iota_{l,l+1}^{\lambda(\Lambda)}(w) \in V_l^{\lambda(\Lambda)}$.
Take 
a word $\mu(w) \in S_{l+1}(\Lambda)$
such that 
$w = [\mu(w)]_{l+1} \in V_{l+1}^{\lambda(\Lambda)}.$
Let 
$v_{\mu(w)}^{l+1} \in V_{l+1}$
be the launching vertex for 
$\mu(w)$ so that
$
\Gamma_{l+1}^-(\mu(w)) = 
\Gamma_{l+1}^-(v_{\mu(w)}^{l+1}).
$ 
Take 
a word
$\mu(w') \in S_l(\Lambda)$
such that
$w' = [\mu(w')]_l$
so that 
$[\mu(w)]_l =[\mu(w')]_l$
and
$
\Gamma_l^-(\mu(w)) = 
\Gamma_l^-(\mu(w')).
$
Hence we have
$
\varPhi_V \circ \iota^{\lambda(\Lambda)}(w) = 
\varPhi_V (w') =  v_{\mu(w')}^l.
$
On the other hand,
one knows that
$\iota\circ \varPhi_V(w) = \iota(v_{\mu(w)}^{l+1}).$
By the local property of
$\lambda$-graph system,
$\mu(w)$ leaves the vertex $\iota(v_{\mu(w)}^{l+1})$.
As
${\frak L}$ is $\lambda$-synchronizing,
$\iota(v_{\mu(w)}^{l+1})$
is the unique vertex which $\mu(w)$ leaves.
Hence
$\iota(v_{\mu(w)}^{l+1})$
is the launching vertex in $V_l$ for
$\mu(w)$.
Since
$
\Gamma_l^-(\mu(w)) = 
\Gamma_l^-(\mu(w')),
$
one sees that 
$\iota(v_{\mu(w)}^{l+1})$
is the launching vertex in $V_l$ for
$\mu(w')$
so that 
$v_{\mu(w')}^l 
= 
\iota(v_{\mu(w)}^{l+1})$.
This means 
$
\varPhi_V \circ \iota^{\lambda(\Lambda)}(w) = 
\iota\circ  
\varPhi_V(w)$.

3.  Injectivity of $\varPhi_V$ :
Suppose that $u, u' \in V_l^{\lambda(\Lambda)}$
satisfy
$\varPhi_V(u) = \varPhi_V(u')$.
Take
$\mu(u), \mu(u') \in S_l(\Lambda)$
such that
$u = [\mu(u)]_l, u' = [\mu(u')]_l$. 
Let
$v_{\mu(u)}^l, v_{\mu(u')}^l $
be the launching vertices in $V_l$ for
$\mu(u), \mu(u')$ 
respectively so that
$
\Gamma_l^-(v_{\mu(u)}^l)
 = 
\Gamma_l^-(\mu(u))$,
$
\Gamma_l^-(v_{\mu(u')}^l) = 
\Gamma_l^-(\mu(u')).
$
Since 
$\varPhi_V(u) = \varPhi_V(u')$,
one sees that 
$v_{\mu(u)}^l = v_{\mu(u')}^l$
so that 
$
\Gamma_l^-(\mu(u))
 = 
\Gamma_l^-(\mu(u'))
$
and
$[\mu(u)]_l = [\mu(u')]_l$.
Hence we have
$u=u'$.

4.  Surjectivity of $\varPhi$ :
For a vertex $v_i^l \in V_l$, 
take a word $\nu \in S_{v_i^l}(\Lambda)$ 
such that $v_i^l$ launches $\nu$.  
By Lemma 3.4,
one sees 
$\nu \in S_l(\Lambda)$.
Hence $\nu$ defines a vertex 
$[\nu]_l$ in
$V_l^{\lambda(\Lambda)}$.
By definition, one has
$\varPhi_V([\nu]_l) = v_i^l$
so that $\varPhi_V$ is surjective.

 5. Existence of edge map $\varPhi_E$ : 
 For $e \in E_{l,l+1}^{\lambda(\Lambda)}$,
put 
$u=s(e) \in V_l^{\lambda(\Lambda)},
v=t(e) \in V_{l+1}^{\lambda(\Lambda)},
\alpha =\lambda(e) \in \Sigma$.
Take a word
$\nu \in S_{l+1}(\Lambda)$
such that $\nu$ starts at $v$.
Put
$\zeta = \alpha \nu \in S_l(\Lambda)$.
As
${\frak L}$ is $\lambda$-synchronizing,
there  exist launching vertices  
$v_{\mu(\zeta)}^l \in V_l$ for $\zeta$
and
$v_{\mu(\nu)}^{l+1} \in V_{l+1}$ for $\nu$
such that
$
\Gamma_l^-(\zeta) = \Gamma_l^-(v_{\mu(\zeta)}^l)
$
and
$
\Gamma_{l+1}^-(\nu) = \Gamma_{l+1}^-(v_{\mu(\nu)}^{l+1})
$
respectively.
Put
$u' = \varPhi_V(u) = v_{\mu(\zeta)}^l \in V_l$
and
$v' = \varPhi_V(v) = v_{\mu(\nu)}^{l+1} \in V_{l+1}$.
Since
$\Gamma_{l+1}^-(\nu) =\Gamma_{l+1}^-(v)$,
one has
$\Gamma_{l+1}^-(v') =\Gamma_{l+1}^-(v)$.
Hence there exist
$e' \in E_{l,l+1}$ in ${\frak L}$
such that
$\lambda(e') = \alpha$
and
$t(e') = v'$.
We set
$\varPhi_E(e) = e'$.
Since ${\frak L}$ is left-resolving,
such $e'$ is unique.
It then follows that
\begin{equation*}
\varPhi_V(t(e))   = v' = t(e') = t(\varPhi_E(e)),
\quad
\varPhi_V(s(e))   = u' = s(e') = s(\varPhi_E(e)).
\end{equation*} 
Both  
${\frak L}^{\lambda(\Lambda)}$ 
and
${\frak L}$
 are left-resolving,
the map $\varPhi_E$ is injective.
Surjectivity of $\varPhi_E$ 
is easily shown so that
 ${\frak L}^{\lambda(\Lambda)}$ and
${\frak L}$
 are isomorphic as $\lambda$-graph systems.
\end{pf}

We call ${\frak L}^{\lambda(\Lambda)}$ 
{\it the}\/ $\lambda$-{\it synchronizing}\/ 
$\lambda$-graph system 
for a $\lambda$-synchronizing subshift $\Lambda$.

\noindent
{\bf Definition.} 
A $\lambda$-graph system
${\frak L}$ is said to be {\it minimal}\/ if
there is no proper $\lambda$-graph subsystem of ${\frak L}$
that presents the same subshift presented by ${\frak L}$.
This means that if 
${\frak L}'$ is a  $\lambda$-graph subsystem of ${\frak L}$
and 
presents the same subshift as the subshift presented by 
${\frak L}$,
then 
${\frak L}'$ coincides with ${\frak L}$.

\begin{prop}
For a $\lambda$-synchronizing subshift $\Lambda$,
the $\lambda$-synchronizing
$\lambda$-graph system 
${\frak L}^{\lambda(\Lambda)}$
is minimal.
\end{prop}
\begin{pf}
Suppose that
${\frak L}' = (V', E', \lambda',\iota')$
be a $\lambda$-graph subsystem of ${\frak L}^{\lambda(\Lambda)}$ 
that presents $\Lambda$.
Hence we have
$
V'_l \subset V_l^{\lambda(\Lambda)}, 
E'_{l,l+1} \subset E_{l,l+1}^{\lambda(\Lambda)}
$
for all
$
 l \in \Zp.
 $
Suppose that 
there exists a vertex $v_i^l \in V_l^{\lambda(\Lambda)}$
such that 
$v_i^l \not\in V'_l$.
By the $\lambda$-synchronization of
${\frak L}^{\lambda(\Lambda)}$,
there exists a synchronizing word
$\mu \in S_l(\Lambda)$
such that
$v_i^l$ launches $\mu$.
There is no any other vertex in $V_l^{\lambda(\Lambda)}$
than $v_i^l$ 
which the word $\mu$ leaves.
Hence the word $\mu$ does not appear in
the presentation of ${\frak L}'$,
 a contradiction.
 Therefore we have
 $V'_l = V_l^{\lambda(\Lambda)}$
 for all $l \in \Zp$.
 We will next show that
$
E'_{l,l+1} = E_{l,l+1}^{\lambda(\Lambda)}
$
for all
$ l \in \Zp$.
Take an arbitrary  edge 
$e \in E_{l,l+1}^{\lambda(\Lambda)}$.
Put
$\alpha =\lambda(e) \in \Sigma$
and
$v_i^l = s(e) \in V_l^{\lambda(\Lambda)}, 
v_j^{l+1} = t(e)\in V_{l+1}^{\lambda(\Lambda)}$. 
Take an $l+1$-synchronizing word 
$\nu \in S_{l+1}(\Lambda)$
such that 
$v_j^{l+1}$ launches $\nu$.
Hence
$\alpha \nu$ leaves $v_i^l$.
Suppose that
$\alpha \nu$ leaves  a vertex $v_{i'}^l$.
Let
$e' \in E_{l,l+1}^{\lambda(\Lambda)}$
be an edge labeled $\alpha$
such that
$s(e') = v_{i'}^l$.
Since $v_j^{l+1}$ is the launching vertex for 
$\nu$,
the word $\nu$ must leave $v_j^{l+1}$
so that 
$v_j^{l+1} = t(e')$.
As ${\frak L}^{\lambda(\Lambda)}$
is left-resolving,
one has 
$e=e'$ and hence
$v_i^l = v_{i'}^l$.
Hence $\alpha \nu$ must leave 
the vertex $v_i^l$.
This means that 
$e$ is the only edge in $E_{l,l+1}$
whose label is the leftmost of
the word $\alpha \nu$.
Therefore 
$e \in E'_{l,l+1}$
and
we have
$E'_{l,l+1} = E_{l,l+1}^{\lambda(\Lambda)}$
for all $l \in \Zp$.
We thus conclude 
${\frak L}' = {\frak L}^{\lambda(\Lambda)}$.
\end{pf}

\section{$\lambda$-synchronization and flow equivalence}
As in Proposition 2.2,
$\lambda$-synchronization is invariant under topological conjugacy.
In this section we will prove that 
$\lambda$-synchronization is invariant even under flow equivalence.
Parry-Sullivan showed that the flow equivalence relation on homeomorphisms
of Cantor sets 
is generated by topological conjugacy and expansion of symbols (\cite{PS}).
Let $\Lambda$ be a subshift over alphabet 
$\Sigma =\{ 1,2,\dots, N \}$.
We define a new subshift $\widetilde{\Lambda}$ 
over the alphabet 
$\widetilde{\Sigma} =\{0, 1,2,\dots, N \}$
as the subshift 
consisting of all biinfinite sequences of $\widetilde{\Sigma}$ obtained by replacing the symbol $1$ in a biinfinite sequence in the subshift $\Lambda$ by the word $01$.  
This operation is called expansion that
corresponds to the equivalence relation called Kakutani equivalence.
An argument in \cite[Proposition]{PS} (cf. \cite[Lemma 2.1]{2001ETDS})
says:
\begin{lem}[\cite{PS}]   
Flow equivalence relation of subshifts is generated by topological conjugacy and the expansion 
$ \Lambda \rightarrow \widetilde{\Lambda}$.
\end{lem}
For a subshift $\Lambda$ over $\Sigma$,
recall that 
$X_{\Lambda} \subset \Sigma^{\Bbb N}$
is defined by the set of all right one-sided sequences
$(x_i)_{i \in {\Bbb N}} \in \Sigma^{\Bbb N}$
such that 
$ 
(x_i)_{i \in {\Bbb Z}} \in \Lambda$,
and 
$ 
X_{\widetilde{\Lambda}} \subset \widetilde{\Sigma}^{\Bbb N}
$
is similarly defined.
We set for $l \in {\Bbb N}$,
\begin{alignat*}{2}
B_{1,l}(\Lambda)  
&=   \{ \mu_1\cdots\mu_l \in B_l(\Lambda)  \mid \mu_1 = 1 \}, &\qquad
 B_{1,*}(\Lambda)  
&= \cup_{l=1}^\infty B_{1,l}(\Lambda),\\
B_{1,l}(\widetilde{\Lambda})  
&=  \{ \nu_1\cdots\nu_l \in B_l(\widetilde{\Lambda}) \mid \nu_1 = 1 \}, &\qquad
B_{1,*}(\widetilde{\Lambda})  
&= \cup_{l=1}^\infty B_{1,l}(\widetilde{\Lambda}),\\
B_{l,0}(\widetilde{\Lambda})  
&=  \{ \nu_1\cdots\nu_l \in B_l(\widetilde{\Lambda}) \mid \nu_l = 0 \},  &\qquad
 B_{*,0}(\widetilde{\Lambda})  
&= \cup_{l=1}^\infty B_{l,0}(\widetilde{\Lambda}),\\
B_{1,l,0}(\widetilde{\Lambda}) 
&= B_{1,l}(\widetilde{\Lambda}) 
\cup B_{l,0}(\widetilde{\Lambda}), &\qquad\quad
B_{1,*,0}(\widetilde{\Lambda})  
&= \cup_{l=1}^\infty B_{1,l,0}(\widetilde{\Lambda}).
\end{alignat*}
Define
$$
\xi^B: B_*(\Lambda)\longrightarrow 
B_*(\widetilde{\Lambda})\backslash B_{1,*,0}(\widetilde{\Lambda})
$$
by  putting the word $01$ in place of $1$ in the words 
$\mu \in B_*(\Lambda)$
from the left in order such as
$$
\xi^B(1 1 2 1 2 1 3 2 1 ) 
= 0 1 0 1 2 0 1 2 0 1 3 2 0 1. 
$$
The symbol $1$ and the symbol $0$ 
can never appear in the leftmost  
and in the rightmost  
of the form $\xi^B(\mu)$ for $\mu \in B_*(\Lambda)$
respectively
so that we have
$$
\xi^B(B_*(\Lambda)) \cap
B_{1,*,0}(\widetilde{\Lambda}) = \emptyset.
$$
We write $\xi^B(\mu)$ as 
$\tilde{\mu}$ for brevity.
Define
$$
\eta^B :
B_*(\widetilde{\Lambda})\backslash B_{1,*,0}(\widetilde{\Lambda}) 
\longrightarrow
B_*(\Lambda)
$$ 
by putting $1$ in place of $01$ in the words 
$\nu \in B_*(\widetilde{\Lambda})\backslash B_{1,*,0}(\widetilde{\Lambda}) 
$
from the left in order such as
$$
\eta^B(0 1 0 1 2 0 1 2 0 1 3 2 0 1)
= 1 1 2 1 2 1 3 2 1.
$$
We write 
$\eta^B(\nu)$ as
$\bar{\nu}$ for brevity.
Hence we have
$$
\eta^B \circ \xi^B = \id_{B_*(\Lambda)},
\qquad 
\xi^B \circ \eta^B = \id_{B_*(\widetilde{\Lambda})\backslash B_{1,*,0}(\widetilde{\Lambda})}.
$$
 In the set 
$S_l(\Lambda)$ of $l$-synchronizing words,
put
\begin{equation*}
S_{1,l}(\Lambda) = \{ \mu_1\cdots\mu_n \in S_l(\Lambda) \mid \mu_1 = 1 \}.
\end{equation*}
We
similarly use the notation 
$S_l(\widetilde{\Lambda})$
for the subshift
$\widetilde{\Lambda}$ 
as the set of $l$-synchronizing words of $\widetilde{\Lambda}$.
Put
\begin{align*}
S_{1,l}(\widetilde{\Lambda}) & = \{ \nu_1\cdots \nu_n \in S_l(\widetilde{\Lambda}) \mid \nu_1 = 1 \},\\
S_{l,0}(\widetilde{\Lambda}) & = \{ \nu_1\cdots \nu_n \in S_l(\widetilde{\Lambda}) \mid \nu_n = 0 \},\\
S_{1,l,0}(\widetilde{\Lambda}) & = S_{1,l}(\widetilde{\Lambda})
\cup S_{l,0}(\widetilde{\Lambda}).
\end{align*}
\begin{lem}
Assume that $\Lambda$ is irreducible.
\begin{enumerate}
\renewcommand{\labelenumi}{(\roman{enumi})}
\item
$\xi^B : B_*(\Lambda) \longrightarrow 
B_*(\widetilde{\Lambda}) \backslash B_{1,*,0}(\widetilde{\Lambda})
$
induces a map
$
\xi^S_l  : S_l(\Lambda) \longrightarrow 
S_l(\widetilde{\Lambda}) \backslash S_{1,l,0}(\widetilde{\Lambda}).
$
\item
$\eta^B :  
B_*(\widetilde{\Lambda}) \backslash B_{1,*,0}(\widetilde{\Lambda})
\longrightarrow B_*(\Lambda) $
induces a map
$
\eta^S_l  :
S_{2l}(\widetilde{\Lambda}) \backslash S_{1,2l,0}(\widetilde{\Lambda}) 
\longrightarrow  S_l(\Lambda).
$
\end{enumerate}
\end{lem}
\begin{pf}
(i) 
Put
$$
X_{\Lambda_1} = \{ (x_1,x_2,\cdots )\in X_\Lambda | x_1 = 1\},
\quad
X_{\widetilde{\Lambda}_1} = \{ (y_1,y_2,\cdots )\in X_{\widetilde{\Lambda}} 
| y_1 = 1\}.
$$
Define
$\xi :X_{\Lambda} \rightarrow 
X_{\widetilde{\Lambda}}\backslash X_{\widetilde{\Lambda}_1}$ 
by putting the word $01$ in place of $1$ in the sequence 
$(x_1,x_2,\cdots ) \in X_{\Lambda}$
from the left in order such as
$$
\xi(1,1,2,1,2,1,3,2,1,1,2,\cdots ) 
=(0,1,0,1,2,0,1,2,0,1,3,2,0,1,0,1,2,\cdots ).
$$
The symbol $1$ can never appear in the first coordinate in the sequences 
of $\xi(X_{\Lambda})$ 
so that we have
$\xi(X_\Lambda) \cap
X_{\widetilde{\Lambda}_1} = \emptyset$.
Define
$
\eta :
X_{\widetilde{\Lambda}}\backslash X_{\widetilde{\Lambda}_1} \rightarrow
X_{\Lambda}
$ 
by putting the word $1$ in place of $01$ in the sequence 
$(y_1,y_2,\cdots ) \in X_{\widetilde{\Lambda}_1}$ with $y_1\ne 1$
from the left in order such as
$$
\eta(0,1,0,1,2,0,1,2,0,1,3,2,0,1,0,1,2,\cdots )
=(1,1,2,1,2,1,3,2,1,1,2,\cdots ) .
$$
Hence we have
$$
\eta \circ \xi = \id_{X_{\Lambda}},
\qquad 
\xi \circ \eta = \id_{X_{\widetilde{\Lambda}}\backslash X_{\widetilde{\Lambda}_1}}.
$$
We will prove that 
if $\mu \in B_*(\Lambda)$ is an $l$-synchronizing word in $\Lambda$, 
the word 
$
\tilde{\mu}(=\xi^B(\mu))
$ 
is an $l$-synchronizing word in $\widetilde{\Lambda}$.
For $y \in X_{\widetilde{\Lambda}}$ with
 $y \in \Gamma_{\infty}^+(\tilde{\mu})$,
 we will show that
\begin{equation*}
 \Gamma_l^-(\tilde{\mu}) \subset \Gamma_l^-(\tilde{\mu}y).
\end{equation*}
Take an arbitrary word 
$\omega=\omega_1\cdots \omega_l \in \Gamma_l^-(\tilde{\mu})$.
 
 Case 1: $\omega_1 \ne 1$.
 
 Since the leftmost of $\omega\tilde{\mu}$
is not $1$
 and  the rightmost of $\omega\tilde{\mu}$
is not $0$,
  replace the symbol $01$ by $1$ in $\omega\tilde{\mu}$,
  one has
  $\eta^B(\omega\tilde{\mu}) \in B_*(\Lambda)$.
Since the leftmost of $\tilde{\mu}$
is not $1$,
the rightmost of 
$\omega$
is not $0$
so that $\omega \not\in B_{1,*,0}(\widetilde{\Lambda})$.
It follows that
$$
\eta^B(\omega \tilde{\mu}) = 
 \eta^B(\omega)\eta^B(\tilde{\mu})= 
 \eta^B(\omega) \mu.
$$
The leftmosts of both $\tilde{\mu}$ and $y$ are not  $1$.
Hence  one has
 $$
 \eta(\tilde{\mu}y) = \eta^B(\tilde{\mu})\eta(y) = \mu \eta(y).
 $$
 Now $\mu$ is $l$-synchronizing in $\Lambda$
 and
 $\eta^B(\omega) \in \Gamma^-_*(\mu)$
with $|\eta^B(\omega)|\le l$.
We have
 $\eta^B(\omega) \mu\eta(y) \in X_\Lambda.$ 
 As
 $\eta^B(\omega)\mu\eta(y)= \eta^B(\omega)\eta^B(\tilde{\mu})\eta(y)$,
 it follows taht
 \begin{equation*}
\omega\tilde{\mu}y = \xi(\eta^B(\omega)\eta^B(\tilde{\mu})\eta(y)) \in 
X_{\widetilde{\Lambda}}\backslash X_{\widetilde{\Lambda}_1}  
 \end{equation*}
 so that 
 $\omega\tilde{\mu}y \in X_{\widetilde{\Lambda}}$
 and hence $\omega \in \Gamma_l^-(\tilde{\mu}y)$.
 
 Case 2: $\omega_1 = 1$.

For $\omega=\omega_1\cdots \omega_l \in \Gamma_l^-(\tilde{\mu})$,
consider
$
\omega'= 0\omega_1\cdots \omega_l \in B_{l+1}(\widetilde{\Lambda}).
$
Since 
$|\eta^B(\omega')|\le l$,
one may apply the above discussion for $\omega'$
so that  
$\omega' \tilde{\mu}y \in X_{\widetilde{\Lambda}}$.
Hence 
$\omega \in \Gamma^-_l(\tilde{\mu}y)$
and
$$
\Gamma^-_l(\tilde{\mu})\subset \Gamma^-_l(\tilde{\mu}y)
\qquad \text{ for all } y \in X_{\widetilde{\Lambda}}
\text{ with } 
y \in \Gamma_{\infty}^+(\tilde{\mu}).
$$
 Therefore 
 $\tilde{\mu}$ is $l$-synchronizing in $\widetilde{\Lambda}$
 and 
$\xi^B: B_*(\Lambda) \longrightarrow
B_*(\widetilde{\Lambda}) \backslash B_{1,*,0}(\widetilde{\Lambda})
$
induces a map
$$
\xi^S_l  : S_l(\Lambda) \longrightarrow 
S_l(\widetilde{\Lambda}) \backslash S_{1,l,0}(\widetilde{\Lambda}) \quad
$$

(ii)
For 
$
\nu \in B_*(\widetilde{\Lambda}) \backslash B_{1,*,0}(\widetilde{\Lambda}),
$
put
$\bar{\nu} =\eta^B(\nu)$.
Suppose that
$\nu$ is $2l$-synchronizing in $\widetilde{\Lambda}$.
We will show that 
$\bar{\nu}$ is $l$-synchronizing in $\Lambda$.
For 
$\gamma \in \Gamma^-_l(\bar{\nu})$,
one sees that
$\xi^B(\gamma\bar{\nu}) \in 
B_*(\widetilde{\Lambda}) \backslash B_{1,*,0}(\widetilde{\Lambda}).
$
Since
$\nu \not\in B_{1,*,0}(\widetilde{\Lambda})$,
one has
$\xi^B(\bar{\nu}) = \nu$
so that
$\xi^B(\gamma)\nu \in B_*(\widetilde{\Lambda})$.
As
$|\gamma| =l$,
one has 
$|\xi^B(\gamma)| \le 2l$.
For 
$x \in \Gamma_\infty^+(\bar{\nu})$,
one has
$$
\xi(\bar{\nu}x) = \nu\xi(x) \in 
X_{\widetilde{\Lambda}}\backslash X_{\widetilde{\Lambda}_1}.
$$
Since $\nu$ is $2l$-synchronizing in $\widetilde{\Lambda}$
and
$\xi^B(\gamma)\nu \in B_*(\widetilde{\Lambda})$
with
$|\xi^B(\gamma)|\le 2l$,
we have
$\xi^B(\gamma)\nu\xi(x) \in X_{\widetilde{\Lambda}}$.
As the leftmost of $\xi^B(\gamma)$ is not $1$,
we have
$\xi^B(\gamma)\nu\xi(x) \in 
X_{\widetilde{\Lambda}}\backslash X_{\widetilde{\Lambda}_1}$.
Hence we have
$$
\gamma \bar{\nu}x
=
\eta(\xi^B(\gamma)\nu\xi(x) ) \in X_\Lambda  
$$
so that
$\gamma \in \Gamma^-_l(\bar{\nu}x)$.
Therefore we have 
$$
\Gamma^-_l(\bar{\nu}) \subset \Gamma^-_l(\bar{\nu}x)
$$
so that
$\bar{\nu}$ is $l$-synchronizing in $\Lambda$.
We thus have 
$\eta^B :  
B_*(\widetilde{\Lambda}) \backslash B_{1,*,0}(\widetilde{\Lambda})
\longrightarrow B_*(\Lambda) 
$
induces 
a map
$$
\eta^S_l  :
S_{2l}(\widetilde{\Lambda}) \backslash S_{1,2l,0}(\widetilde{\Lambda}) 
\longrightarrow  S_l(\Lambda).
$$
\end{pf}
We note  that an irreducible subshift $\Lambda$ is $\lambda$-synchronizing
if and only if for 
$
\mu\in B_*(\Lambda)
$ 
and
$
k \in {\Bbb N}
$
there exists
$\nu \in S_k(\Lambda)
$
such that
$
\mu \nu \in B_*(\Lambda).
$
\begin{prop}
An irreducible subshift
$\Lambda$
is $\lambda$-synchronizing if and only if 
so is $\widetilde{\Lambda}$.
\end{prop}
\begin{pf}
It is easy to see that
$\Lambda$ is irreducible if and only if 
so is $\widetilde{\Lambda}$.
Suppose that 
$\Lambda$ is $\lambda$-synchronizing.
For $\mu=\mu_1\cdots \mu_l \in B_*(\widetilde{\Lambda})$
and
$k \in {\Bbb N}$,
we have three cases.

Case 1: $\mu_1 \ne 1, \, \mu_l \ne 0$.

As 
$\mu \in B_*(\widetilde{\Lambda})\backslash B_{1,*,0}(\widetilde{\Lambda})$,
by putting $1$ in place of $01$,
we have 
$\bar{\mu} =\eta^B(\mu) \in B_*(\Lambda)$.
Since $\Lambda$ is $\lambda$-synchronizing,
there exists $\nu \in S_k(\Lambda)$
such that
$\bar{\mu}\nu \in B_*(\Lambda)$.
Put
$\tilde{\nu} =\xi^B(\nu) \in S_k(\widetilde{\Lambda})$
so that one sees 
$$
\mu \tilde{\nu} 
=\xi^B(\bar{\mu})\xi^B(\nu)=\xi^B(\bar{\mu}\nu) \in B_*(\widetilde{\Lambda}).
$$

Case 2: $\mu_1 =1, \, \mu_l \ne 0$.

Consider the word
$0\mu = 0\mu_1\cdots \mu_l \in B_*(\widetilde{\Lambda})$
so that 
$
0\mu \in B_*(\widetilde{\Lambda})\backslash B_{1,*,0}(\widetilde{\Lambda}).
$
By the above discussion of Case 1 for the word $0\mu$,
there exists $\nu \in S_k(\Lambda)$
such that 
$0\mu \xi^B(\nu) \in B_*(\widetilde{\Lambda})$.
Put
$\tilde{\nu} = \xi^B(\nu)$
so that 
$\mu \tilde{\nu} \in B_*(\widetilde{\Lambda})$.

Case 3: $\mu_l =0$.

Put
$\mu 1 = \mu_1 \cdots \mu_l 1 \in B_*(\widetilde{\Lambda})$.
As
the rightmost of $\mu 1$ is not $0$,
by applying the above two cases to $\mu 1$ and $k+1$,
one finds 
$\tilde{\nu} \in S_{k+1}(\widetilde{\Lambda})$
such that
$\mu 1 \tilde{\nu} \in B_*(\widetilde{\Lambda})$.
Put
$\hat{\nu} = 1 \tilde{\nu} \in S_k(\widetilde{\Lambda})$
so that we have
$
\mu \hat{\nu} \in B_*(\widetilde{\Lambda}).
$

Therefore we conclude that 
$\widetilde{\Lambda}$ is $\lambda$-synchronizing.

Conversely 
assume that $\widetilde{\Lambda}$ is $\lambda$-synchronizing.
For $\mu=\mu_1\cdots \mu_l \in B_*(\Lambda)$
and
$k \in {\Bbb N}$,
put
$\tilde{\mu} = \xi^B(\mu) \in 
B_*(\widetilde{\Lambda})\backslash B_{1,*,0}(\widetilde{\Lambda})
$.
By $\lambda$-synchronization of $\widetilde{\Lambda}$,
for $\tilde{\mu}$ and $2 k\in {\Bbb N}$,
there exists
$\nu'=\nu'_1\cdots \nu'_n \in S_{2k}(\widetilde{\Lambda})$
such that
$\tilde{\mu}\nu' \in B_*(\widetilde{\Lambda})$.
If $\nu'_n =0$,
consider the word 
$\nu' 1 = \nu'_1\cdots \nu'_n 1$
instead of $\nu'$
so that 
we may assume that 
$\nu'_n \ne 0$.
As $\tilde{\mu}$ does not end at $0$,
$\nu'$ does not begin with $1$
so that
$\nu' \not\in B_{1,*,0}(\widetilde{\Lambda})$.
As
both
$
\tilde{\mu}, \nu' \in 
B_{*}(\widetilde{\Lambda})\backslash B_{1,*,0}(\widetilde{\Lambda})
$,
it follows that
$$
\mu \eta^B(\nu')
=
\eta^B(\tilde{\mu} \nu') \in B_*(\Lambda).
$$
By putting 
$\nu = \eta^B(\nu') $,
one sees that
$\nu \in S_k(\Lambda)$
and
$\mu \nu \in B_*(\Lambda)$.
This shows that  
$\Lambda$ is $\lambda$-synchronizing.
\end{pf}
Therefore we conclude 
\begin{thm}
The $\lambda$-synchronization is invariant under flow equivalence 
of subshifts.
\end{thm}
\begin{pf}
By Proposition 2.2 and Proposition 4.3,
$\lambda$-synchronization is invariant under topological conjugacy and expansions of subshifts.
Therefore by Lemma 4.1,
$\lambda$-synchronization is invariant under flow equivalence 
of subshifts.
\end{pf}

\section{K-groups and Bowen-Franks groups}

In this section, we will prove that 
the K-groups and the Bowen-Franks groups 
for the $\lambda$-synchronizing $\lambda$-graph system for $\lambda$-synchronizing subshifts
are invariant under expansion
 $\Lambda \longrightarrow \widetilde{\Lambda}$. 
The line of the proof basically follows the proof of
\cite[Theorem]{2001ETDS}.
As a result, the groups yield invariants for 
flow equivalence of $\lambda$-synchronizing subshifts.
The K-groups and the Bowen-Franks groups for
$\lambda$-synchronizing $\lambda$-graph systems
are defined to be those groups 
for the $\lambda$-graph system
$\LLL$ ( \cite{KM2010}).
They are called the $\lambda$-synchronizing K-groups for $\Lambda$
and
the $\lambda$-synchronizing Bowen-Franks groups for $\Lambda$
respectively.
We will briefly describe them.
Let
$m_\lambda(l)$ be the cardinal number of
the vertex set
$V_l^{\lambda(\Lambda)}$
denoted by 
$\{ v_1^l, \dots,v_{m_{\lambda}(l)}^l \}$.
Define $m_\lambda(l) \times m_\lambda(l+1)$ matrices
 $I^{\lambda(\Lambda)}_{l,l+1}$ and   
$A^{\lambda(\Lambda)}_{l,l+1}$
by setting
\begin{align*}
I^{\lambda(\Lambda)}_{l,l+1}(i,j)
& =
{
\begin{cases}
1 & \text{ if } \iota^{\lambda(\Lambda)}(v_j^{l+1}) = v_i^l, \\
0 & \text{ otherwise, }
\end{cases}
} \\
A^{\lambda(\Lambda)}_{l,l+1}(i,j)
& =
\text{ the number of the labeled edges from } v_i^l \text{ to } v_j^{l+1}
\end{align*}
for
$i=1,\dots,m_\lambda(l)$ and 
$j=1,\dots,m_\lambda(l+1)$.
The sequence
$(A^{\lambda(\Lambda)}_{l,l+1},I^{\lambda(\Lambda)}_{l,l+1}), l \in \Zp$
of pairs of matrices becomes a nonnegative matrix system (\cite{1999DocMath}).
The $\lambda$-synchronizing K-groups $K_i^\lambda(\Lambda), i=0,1$
are defined 
as the K-groups 
$K_i(A^{\lambda(\Lambda)}, I^{\lambda(\Lambda)}), i=0,1$
for the nonnegative matrix system $(A^{\lambda(\Lambda)}, I^{\lambda(\Lambda)})$
that are formulated by
\begin{align*}
K_0^{\lambda}(\Lambda) 
& = 
\underset{l}{\varinjlim} \{
{}^t\negthinspace \bar{I}^{\lambda(\Lambda)}_{l,l+1}
:
\Coker({}^t\negthinspace I^{\lambda(\Lambda)}_{l-1,l} -
{}^t\negthinspace A^{\lambda(\Lambda)}_{l-1,l} )
\longrightarrow
\Coker({}^t\negthinspace I^{\lambda(\Lambda)}_{l,l+1} -
{}^t\negthinspace A^{\lambda(\Lambda)}_{l,l+1} )
\}, \\
K_1^{\lambda}(\Lambda) 
& = 
\underset{l}{\varinjlim}
\{ {}^t\negthinspace I^{\lambda(\Lambda)}_{l,l+1} :
\Ker  
({}^t\negthinspace I^{\lambda(\Lambda)}_{l-1,l} - 
{}^t\negthinspace A^{\lambda(\Lambda)}_{l-1,l} )
\longrightarrow
\Ker  
({}^t\negthinspace I^{\lambda(\Lambda)}_{l,l+1} - 
{}^t\negthinspace A^{\lambda(\Lambda)}_{l,l+1} ) \},
\end{align*}
where
\begin{align*}
\Coker({}^t\negthinspace I^{\lambda(\Lambda)}_{l,l+1} -
{}^t\negthinspace A^{\lambda(\Lambda)}_{l,l+1} ) 
& 
=
{\Bbb Z}^{m_\lambda(l+1)} /  
({}^t\negthinspace I^{\lambda(\Lambda)}_{l,l+1} -
{}^t\negthinspace A^{\lambda(\Lambda)}_{l,l+1} ){\Bbb Z}^{m_\lambda(l)}, \\
\Ker  
({}^t\negthinspace I^{\lambda(\Lambda)}_{l,l+1} - 
{}^t\negthinspace A^{\lambda(\Lambda)}_{l,l+1} )
&
=
\Ker  
({}^t\negthinspace I^{\lambda(\Lambda)}_{l,l+1} - 
{}^t\negthinspace A^{\lambda(\Lambda)}_{l,l+1} )
\quad \text{ in } {\Bbb Z}^{m_\lambda(l)},
\end{align*}
and
${}^t\negthinspace\bar{I}^{\lambda(\Lambda)}_{l,l+1}
$
is the natural homomorphism induced by 
${}^t\negthinspace I^{\lambda(\Lambda)}_{l,l+1}: {\Bbb Z}^{m_\lambda(l)} \longrightarrow 
{\Bbb Z}^{m_\lambda(l+1)}.
$
Denote by
$
{\Bbb Z}^{\lambda(\Lambda)}
$
the projective limit  
$
\varprojlim \{I^{\lambda(\Lambda)}_{l,l+1}: {\Bbb Z}^{m_\lambda(l+1)} \rightarrow 
                                      {\Bbb Z}^{m_\lambda(l)} \}
$ 
of abelian group.
The sequence $A^{\lambda(\Lambda)}_{l,l+1}, l \in \Zp$ 
of matrices acts on ${\Bbb Z}^{\lambda(\Lambda)}$
as an endomorphism that we denote by
$A^{\lambda(\Lambda)}$.
The identity on ${\Bbb Z}^{\lambda(\Lambda)}$ 
is denoted by $I$.
The $\lambda$-synchronizing Bowen-Franks groups 
are formulated by  
$$
BF_{\lambda(\Lambda)}^0(\Lambda) 
= {\Bbb Z}^{\lambda(\Lambda)} /(I-A^{\lambda(\Lambda)}){\Bbb Z}^{\lambda(\Lambda)}, 
\qquad
BF_{\lambda(\Lambda)}^1(\Lambda) 
= \Ker(I-A^{\lambda(\Lambda)}) \text{ in }{\Bbb Z}^{\lambda(\Lambda)}.
$$
The group
$BF_{\lambda(\Lambda)}^0(\Lambda) $
is computed from the K-groups 
$K_i^\lambda(\Lambda), i=0,1$ by the 
 universal coefficient  type theorem (\cite[Theorem 9.6]{1999DocMath})
described as 
$$
0
\longrightarrow
\Ext_{\Bbb Z}^1(K^{\lambda(\Lambda)}_0(\Lambda),{\Bbb Z})
\overset{\delta}{\longrightarrow}
BF_{\lambda(\Lambda)}^0(\Lambda) 
\overset{\gamma}{\longrightarrow}
 \Hom_{\Bbb Z}(K^{\lambda(\Lambda)}_1(\Lambda),{\Bbb Z})
\longrightarrow
0
$$
that splits unnaturally.
For the group 
$BF_{\lambda(\Lambda)}^1(\Lambda)$, 
there exists an isomorphism (\cite{1999DocMath}):  
$$
BF_{\lambda(\Lambda)}^1(\Lambda) 
\cong \Hom_{\Bbb Z}(K^{\lambda(\Lambda)}_0(\Lambda),{\Bbb Z}).
$$

We henceforth fix a $\lambda$-synchronizing subshift $\Lambda$ over
alphabet $\Sigma =\{1,2,\dots,N\}$.
Recall that 
$S_l(\Lambda)$ denotes the set of $l$-synchronizing words,
and
$S_{1,l}(\Lambda)$ denotes
the subset 
$\{ \mu_1\cdots\mu_n \in S_l(\Lambda) \mid \mu_1 = 1 \}
$
of $S_l(\Lambda)$.
Set
\begin{equation*}
\Omega_\lambda^l(\Lambda)  = S_l(\Lambda) / \sim_l 
\quad \text{ and } \quad
\Omega_\lambda^l(\Lambda_1)  = S_{1,l}(\Lambda) / \sim_l   
\end{equation*}
the $l$-past equivalence classes of
$S_l(\Lambda)$ 
and
$S_{1,l}(\Lambda)$
respectively.  
Similarly for the subshift $\widetilde{\Lambda}$,
set
\begin{equation*}
\Omega_\lambda^l(\widetilde{\Lambda})  = S_l(\widetilde{\Lambda}) / \sim_l 
\quad \text{ and } \quad 
\Omega_\lambda^l(\widetilde{\Lambda}_1)  = S_{1,l}(\widetilde{\Lambda}) / \sim_l
\end{equation*}   
the $l$-past equivalence classes of $S_l(\widetilde{\Lambda})$ 
and
$S_{1,l}(\widetilde{\Lambda})$ 
respectively.
Denote by
$(S_l(\widetilde{\Lambda}) \backslash S_{1,l,0}(\widetilde{\Lambda}) )/ \sim_l 
$
the $l$-past equivalence classes of  
$S_l(\widetilde{\Lambda}) \backslash S_{1,l,0}(\widetilde{\Lambda})$.
\begin{lem}
$
(S_l(\widetilde{\Lambda}) \backslash S_{1,l,0}(\widetilde{\Lambda}) )/ \sim_l 
=
\Omega_\lambda^l(\widetilde{\Lambda}) \backslash \Omega_\lambda^l(\widetilde{\Lambda}_1).
$
\end{lem}
\begin{pf}
Since
$S_l(\widetilde{\Lambda}) \backslash S_{1,l,0}(\widetilde{\Lambda}) 
\subset
S_l(\widetilde{\Lambda})$,
one has
\begin{equation*}
(S_l(\widetilde{\Lambda}) \backslash S_{1,l,0}(\widetilde{\Lambda}) )/ \sim_l 
\subset
\Omega_\lambda^l(\widetilde{\Lambda}).
\end{equation*}
For $\mu = \mu_1\cdots\mu_n \in 
S_l(\widetilde{\Lambda}) \backslash S_{1,l,0}(\widetilde{\Lambda})$,
as
$\mu_1 \ne 1$
one has 
$\nu \not\in\Gamma_l^-(\mu)$
for
all $\nu \in B_{l,0}(\widetilde{\Lambda})$.
Since
for $\eta \in S_{1,l}(\widetilde{\Lambda})$,
one has 
$\nu \in \Gamma_l^-(\eta)$
for some
$\nu \in B_{l,0}(\widetilde{\Lambda})$.
Hence
$[\mu]_l \not\in \Omega_\lambda^l(\widetilde{\Lambda}_1)$
so that
\begin{equation*}
(S_l(\widetilde{\Lambda}) \backslash S_{1,l,0}(\widetilde{\Lambda}) )/ \sim_l 
\subset
\Omega_\lambda^l(\widetilde{\Lambda}) \backslash \Omega_\lambda^l(\widetilde{\Lambda}_1).
\end{equation*}
On the other hand,
for
$[\mu]_l \in  \Omega_\lambda^l(\widetilde{\Lambda}) \backslash \Omega_\lambda^l(\widetilde{\Lambda}_1)$
with $\mu = \mu_1\cdots\mu_n$,  
one has
$\mu_1 \ne 1$.
If $\mu_n =0$,
consider
$\mu 1 = \mu_1\cdots\mu_n 1$
so that
$[\mu]_l = [\mu 1]_l$
and
$\mu 1 \not\in S_{1,l,0}(\widetilde{\Lambda}).$
Hence
$[\mu 1]_l \in (S_l(\widetilde{\Lambda}) \backslash S_{1,l,0}(\widetilde{\Lambda}) )/ \sim_l.$ 
Therefore we have
\begin{equation*}
(S_l(\widetilde{\Lambda}) \backslash S_{1,l,0}(\widetilde{\Lambda}) )/ \sim_l 
=
\Omega_\lambda^l(\widetilde{\Lambda}) \backslash \Omega_\lambda^l(\widetilde{\Lambda}_1).
\end{equation*}
\end{pf}

\begin{lem}
Assume that $\Lambda$ is irreducible.
\begin{enumerate}
\renewcommand{\labelenumi}{(\roman{enumi})}
\item
$
\xi^S_l  : S_l(\Lambda) \longrightarrow 
S_l(\widetilde{\Lambda}) \backslash S_{1,l,0}(\widetilde{\Lambda})
$
induces a map
$
 \xi^\lambda_l  : \Omega_\lambda^l(\Lambda) \longrightarrow 
\Omega_\lambda^l(\widetilde{\Lambda}) \backslash \Omega_\lambda^l(\widetilde{\Lambda}_1). 
$
\item
$
\eta^S_l  :
S_{2l}(\widetilde{\Lambda}) \backslash S_{1,2l,0}(\widetilde{\Lambda}) 
\longrightarrow  S_l(\Lambda)
$
induces a map
$ \eta^\lambda_l  : 
\Omega_\lambda^{2l}(\widetilde{\Lambda}) \backslash \Omega_\lambda^{2l}(\widetilde{\Lambda}_1) \longrightarrow  \Omega_\lambda^l(\Lambda). 
$
\end{enumerate}
\end{lem}
\begin{pf}
(i) 
We will  prove that
for $\mu, \nu \in B_*(\Lambda)$ 
$\mu \sim_{l} \nu$ in $\Lambda$ implies 
$\tilde{\mu} \sim_{l} \tilde{\nu}$ in 
$\widetilde{\Lambda}$.
Suppose that  
 $\mu \sim_{l} \nu$ in $\Lambda$ .
Take an arbitrary word 
$\omega=\omega_1\cdots \omega_l \in \Gamma_l^-(\tilde{\mu})$.

 Case 1: $\omega_1 \ne 1$.
 
  Since the leftmost of $\omega\tilde{\mu}$
is not $1$ 
and the rightmost of $\omega\tilde{\mu}$
is not $0$,
 replace the symbol $01$ by $1$ in $\omega\tilde{\mu}$,
  one has
  $\eta^B(\omega\tilde{\mu})=\eta^B(\omega)\mu\in B_*(\Lambda)$.
As
$|\eta^B(\omega)| \le l$
and
$\Gamma^-_l(\mu) = \Gamma^-_l(\nu)$,
we have
$\eta^B(\omega) \in \Gamma^-_*(\nu)$.
Hence
$\xi^B(  \eta^B(\omega)\nu) \in B_*(\widetilde{\Lambda})$.
Since the leftmost of $\tilde{\mu}$ is not $1$,
the rightmost of $\omega$ is not $ 0$
so that $\omega \not\in B_{1,*,0}(\widetilde{\Lambda})$
and
$
\xi^B(\eta^B(\omega) \nu) = \omega \tilde{\nu}.
$
Hence we have
$\omega \in \Gamma_l^-(\tilde{\nu})$.

 Case 2: $\omega_1 = 1$.

For $\omega=\omega_1\cdots \omega_l \in \Gamma_l^-(\tilde{\mu})$,
consider
$
\omega'= 0\omega_1\cdots \omega_l \in B_{l+1}(\widetilde{\Lambda}).
$
Since 
$|\eta^B(\omega')|\le l$,
one may apply the above discussion for $\eta^B(\omega')$
so that 
the condition
$\eta^B(\omega')\mu \in B_*(\Lambda)$
implies
$\eta^B(\omega')\nu \in B_*(\Lambda)$
and
$\xi^B(\eta^B(\omega')\nu) \in B_*(\widetilde{\Lambda})$.
Hence we have
$\omega'\tilde{\nu} \in B_*(\widetilde{\Lambda})$
so that
$\omega\tilde{\nu} \in B_*(\widetilde{\Lambda})$
and
$\omega \in \Gamma^-_l(\tilde{\nu})$.
We thus conclude 
$\Gamma^-_l(\tilde{\mu})\subset \Gamma^-_l(\tilde{\nu})$
and
then
$\Gamma^-_l(\tilde{\mu}) = \Gamma^-_l(\tilde{\nu})$.

Therefore 
by the preceding lemma,
the map
$\xi_l^S: S_l(\Lambda) 
\longrightarrow
S_l(\widetilde{\Lambda}) \backslash S_{1,l,0}(\widetilde{\Lambda})
$
induces a map
$ \xi^\lambda_l  : \Omega_\lambda^l(\Lambda) \longrightarrow 
\Omega_\lambda^l(\widetilde{\Lambda}) \backslash \Omega_\lambda^l(\widetilde{\Lambda}_1)
$
defined by 
$\xi^\lambda_l([\mu]_l) = [\tilde{\mu}]_l.$

(ii)
It is easy to seee that 
for 
$
\mu,\nu \in 
S_{2l}(\widetilde{\Lambda}) \backslash S_{1,2l,0}(\widetilde{\Lambda})
$
the condition 
$\Gamma^-_{2l}(\mu) = \Gamma^-_{2l}(\nu)$ in $\widetilde{\Lambda}$
implies 
$\Gamma^-_{l}(\bar{\mu}) = \Gamma^-_{l}(\bar{\nu})$ in $\Lambda$.
Hence by the preceding lemma, 
the map
$
\eta^S_l  :
S_{2l}(\widetilde{\Lambda}) \backslash S_{1,2l,0}(\widetilde{\Lambda}) 
\longrightarrow  S_l(\Lambda) 
$
induces a map
$
 \eta^\lambda_l  : 
\Omega_\lambda^{2l}(\widetilde{\Lambda}) \backslash \Omega_\lambda^{2l}(\widetilde{\Lambda}_1) \longrightarrow  \Omega_\lambda^l(\Lambda) 
$
defined by 
$\eta^\lambda_l([\mu]_{2l}) = [\bar{\mu}]_l.$
\end{pf}

The restriction of the natural surjection
$$
\iota^{\lambda(\Lambda)}_{l,l+1}:
[\mu]_{l+1} \in \Omega^{l+1}_\lambda(\Lambda)
\longrightarrow
[\mu]_l \in \Omega^l_\lambda(\Lambda)
$$
yields the surjection
$$
\Omega^{l+1}_\lambda(\Lambda_1)
\longrightarrow
\Omega^l_\lambda(\Lambda_1),
$$
which we still denote by
$\iota^{\lambda(\Lambda)}_{l,l+1}$.
Similarly 
we have the natural surjections 
$$
\Omega^{l+1}_\lambda(\widetilde{\Lambda}_1)
\longrightarrow
\Omega^l_\lambda(\widetilde{\Lambda}_1),
\qquad
\Omega^{l+1}_\lambda(\widetilde{\Lambda})
\backslash\Omega^{l+1}_\lambda(\widetilde{\Lambda}_1)
\longrightarrow
\Omega^l_\lambda(\widetilde{\Lambda})
\backslash\Omega^l_\lambda(\widetilde{\Lambda}_1)
$$
by restricting  
the surjection 
$$
\iota^{\lambda(\widetilde{\Lambda})}_{l,l+1}:
\Omega^{l+1}_\lambda(\widetilde{\Lambda})
\longrightarrow
 \Omega^l_\lambda(\widetilde{\Lambda})
$$
which we still denote by
$\iota^{\lambda(\widetilde{\Lambda})}_{l,l+1}$.
We set the projective limits of the compact Hausdorff spaces: 
\begin{align*}
\Omega_\lambda(\Lambda)
= &
\underset{\iota^{\lambda(\Lambda)}_{l,l+1}}{\varprojlim}\Omega^{l}_\lambda(\Lambda),
\qquad
\Omega_\lambda(\widetilde{\Lambda})
= 
\underset{\iota^{\lambda(\widetilde{\Lambda})}_{l,l+1}}{\varprojlim}\Omega^{l}_\lambda(\widetilde{\Lambda}),\\
\Omega_\lambda(\Lambda_1)
= &
\underset{\iota^{\lambda(\Lambda)}_{l,l+1}}{\varprojlim}\Omega^{l}_\lambda(\Lambda_1),
\qquad
\Omega_\lambda(\widetilde{\Lambda}_1)
= 
\underset{\iota^{\lambda(\widetilde{\Lambda})}_{l,l+1}}{\varprojlim}\Omega^{l}_\lambda(\widetilde{\Lambda}_1).
\end{align*}
We note  
\begin{equation*}
\underset{\iota^{\lambda(\widetilde{\Lambda})}_{l,l+1}}{\varprojlim}
\Omega^l_{\lambda}(\widetilde{\Lambda})\backslash\Omega^l_{\lambda}(\widetilde{\Lambda}_1)
\cong
\underset{\iota^{\lambda(\widetilde{\Lambda})}_{l,l+1}}{\varprojlim}
\Omega^l_{\lambda}(\widetilde{\Lambda})
\backslash
\underset{\iota^{\lambda(\widetilde{\Lambda})}_{l,l+1}}{\varprojlim}
\Omega^l_{\lambda}(\widetilde{\Lambda}_1)
=\Omega_{\lambda}(\widetilde{\Lambda}) \backslash \Omega_{\lambda}(\widetilde{\Lambda}_1).
\end{equation*}
We put
$$
\iota^{\lambda(\Lambda)}_{l,l+n}
=\iota^{\lambda(\Lambda)}_{l,l+1} \circ \iota^{\lambda(\Lambda)}_{l+1,l+2}
\circ \cdots \circ \iota^{\lambda(\Lambda)}_{l+n-1,l+n}:
\Omega^{l+n}_\lambda(\Lambda)
\longrightarrow
\Omega^l_\lambda(\Lambda)
$$
and
$
\iota^{\lambda(\widetilde{\Lambda})}_{l,l+n}:
\Omega^{l+n}_\lambda(\widetilde{\Lambda})
\longrightarrow
 \Omega^l_\lambda(\widetilde{\Lambda})
$
is similarly defined.
The following lemma is straightforward.
\begin{lem}
\hspace{3cm}
\begin{enumerate}
\renewcommand{\labelenumi}{(\roman{enumi})}
\item
$
\eta^\lambda_l\circ\xi^\lambda_{2l}
=\iota^{\lambda(\Lambda)}_{l,2l}: 
\Omega^{2l}_\lambda(\Lambda)
\longrightarrow 
\Omega^{l}_\lambda(\Lambda).
$
\item
$\xi^\lambda_l\circ \eta^\lambda_{l}
=\iota^{\lambda(\widetilde{\Lambda})}_{l,2l}: 
\Omega^{2l}_{\lambda}(\widetilde{\Lambda})\backslash \Omega^{2l}_{\lambda}(\widetilde{\Lambda}_1)
\longrightarrow 
\Omega^l_{\lambda}(\widetilde{\Lambda})\backslash
\Omega^l_{\lambda}(\widetilde{\Lambda}_1).
$
\end{enumerate}
\end{lem}

\begin{cor}
The sequence
$\{ \xi_l^\lambda \}_{l \in {\Bbb N}}$
of the maps
$
\xi_l^\lambda:
\Omega^{l}_\lambda(\Lambda)
\longrightarrow
\Omega^l_{\lambda}(\widetilde{\Lambda})\backslash\Omega^l_{\lambda}(\widetilde{\Lambda}_1),
l \in {\Bbb N}
$
induces a homeomorphism
\begin{equation*}
\underset{l}{\varprojlim}\xi_l^\lambda:
\Omega_\lambda(\Lambda)
=
\underset{\iota^{\lambda(\Lambda)}_{l,l+1}}{\varprojlim}\Omega^{l}_\lambda(\Lambda)
\longrightarrow
\Omega_{\lambda}(\widetilde{\Lambda})\backslash\Omega_{\lambda}(\widetilde{\Lambda}_1)
=
\underset{\iota^{\lambda(\widetilde{\Lambda})}_{l,l+1}}{\varprojlim}
\Omega^l_{\lambda}(\widetilde{\Lambda})\backslash\Omega^l_{\lambda}(\widetilde{\Lambda}_1)
\end{equation*}
of the projective limits as compact Hausdorff spaces.
\end{cor}
We denote by
$\varPhi_d$
the homeomorphism
$
\underset{l}{\varprojlim}\xi_l^\lambda:
\Omega_\lambda(\Lambda)
\longrightarrow
\Omega_{\lambda}(\widetilde{\Lambda})\backslash\Omega_{\lambda}(\widetilde{\Lambda}_1).
$

Let 
$\varphi^B:
B_{1,*}(\Lambda)\longrightarrow B_{1,*}(\widetilde{\Lambda})\backslash B_{*,0}(\widetilde{\Lambda})
$
be
the map defined by putting $01$ in place of $1$
in the subword $x_2\cdots 
x_n$ 
of a word 
$1 x_2\cdots x_n \in B_{1,*}(\Lambda)$
such as 
\begin{equation*}
\varphi^B(1 1 2 1 3 2 1 3 1)
     = 1 0 1 2 0 1 3 2 0 1 3 0 1. 
\end{equation*}
Let 
$\psi^B:
B_{1,*}(\widetilde{\Lambda})\backslash B_{*,0}(\widetilde{\Lambda})
\longrightarrow B_{1,*}(\Lambda)
$
be
the map defined by putting $1$ in place of $01$
in the subword $y_2\cdots y_k$ in a word 
$1 y_2\cdots y_k \in B_{1,*}(\widetilde{\Lambda})$
such as 
\begin{equation*}
\psi^B (1 0 1 2 0 1 3 2 0 1 3 0 1)
   = 1 1 2 1 3 2 1 3 1. 
\end{equation*}

The following lemmas are similarly shown to Lemma 4.2 and Lemma 5.2.
\begin{lem}
\hspace{3cm}
\begin{enumerate}
\renewcommand{\labelenumi}{(\roman{enumi})}
\item
$\varphi^B:
B_{1,*}(\Lambda)
\longrightarrow 
B_{1,*}(\widetilde{\Lambda})\backslash B_{*,0}(\widetilde{\Lambda})
$
induces the maps
\begin{equation*}
\varphi^S_l: S_{1,l-1}(\Lambda)
\longrightarrow S_{1,l}(\widetilde{\Lambda})
\quad
\text{ and }
\quad
\varphi^\lambda_l: \Omega^{l-1}_\lambda(\Lambda_1)
\longrightarrow 
\Omega^l_{\lambda}(\widetilde{\Lambda}_1).
\end{equation*}
\item
$\psi^B:
B_{1,*}(\widetilde{\Lambda})\backslash B_{*,0}(\widetilde{\Lambda})
\longrightarrow B_{1,*}(\Lambda)
$
induces the maps
\begin{equation*}
\psi^S_l: S_{1,2l}(\widetilde{\Lambda})
\longrightarrow 
S_{1,l}(\Lambda)
\quad
\text{ and }
\quad
\psi^\lambda_l: \Omega^{2l}_{\lambda}(\widetilde{\Lambda}_1)
\longrightarrow \Omega^{l}_\lambda(\Lambda_1).
\end{equation*}
\end{enumerate}
\end{lem}

\begin{lem}
\hspace{3cm}
\begin{enumerate}
\renewcommand{\labelenumi}{(\roman{enumi})}
\item
$
\psi^\lambda_l\circ\varphi^\lambda_{2l}
=\iota^{\lambda(\Lambda)}_{2l-1,l}: 
\Omega^{2l-1}_\lambda(\Lambda_1)
\longrightarrow 
\Omega^{l}_\lambda(\Lambda_1).
$
\item
$\varphi^\lambda_l\circ \psi^\lambda_{l-1}
=\iota^{\lambda(\widetilde{\Lambda})}_{2(l-1),l}: 
\Omega^{2(l-1)}_{\lambda}(\widetilde{\Lambda}_1)
\longrightarrow 
\Omega^l_{\lambda}(\widetilde{\Lambda}_1).
$
\end{enumerate}
\end{lem}

\begin{cor}
The sequence 
$\{ \varphi_l^\lambda \}_{l \in \Zp}$
of the maps
$
\varphi_l^\lambda :
\Omega^{l-1}_\lambda(\Lambda_1)
\longrightarrow
\Omega^l_{\lambda}(\widetilde{\Lambda}_1), l \in \Zp$
induces a homeomorphism
\begin{equation*}
\underset{l}{\varprojlim}\varphi_l^\lambda :
\Omega_\lambda(\Lambda_1)
= \underset{\iota^{\lambda(\Lambda)}_{l,l+1}}{\varprojlim}\Omega^{l}_\lambda(\Lambda_1)
\longrightarrow
\Omega_{\lambda}(\widetilde{\Lambda}_1)
=
\underset{\iota^{\lambda(\widetilde{\Lambda})}_{l,l+1}}{\varprojlim}\Omega^l_{\lambda}(\widetilde{\Lambda}_1)
\end{equation*}
of the projective limits as compact Hausdorff spaces.
\end{cor}
We denote by 
$\varPhi_1
$
the homeomorphism
$
\underset{l}{\varprojlim}\varphi_l^\lambda:
\Omega_\lambda(\Lambda_1)
\longrightarrow
\Omega_{\lambda}(\widetilde{\Lambda}_1).
$
By Corollary 5.4 and Corollary 5.7,
we have homeomorphisms
\begin{equation*}
\varPhi_d : 
\Omega_\lambda(\Lambda) 
 \longrightarrow
\Omega_{\lambda}(\widetilde{\Lambda}) \backslash \Omega_{\lambda}(\widetilde{\Lambda}_1),
\qquad
\varPhi_1 :
\Omega_\lambda(\Lambda_1) 
 \longrightarrow
 \Omega_{\lambda}(\widetilde{\Lambda}_1).
\end{equation*}

Therefore we have
\begin{prop}
The disjoint union
$\Omega_\lambda(\Lambda_1) 
\sqcup
\Omega_\lambda(\Lambda)
$
is homeomorphic to
$\Omega_{\lambda}(\widetilde{\Lambda}) $
through 
the homeomorphism
\begin{equation*}
 \varPhi_1 \sqcup\varPhi_d :
\Omega_\lambda(\Lambda_1) 
\sqcup
\Omega_\lambda(\Lambda)
\longrightarrow 
\Omega_{\lambda}(\widetilde{\Lambda}_1)
\sqcup
(\Omega_{\lambda}(\widetilde{\Lambda}) \backslash \Omega_{\lambda}(\widetilde{\Lambda}_1))
=
\Omega_{\lambda}(\widetilde{\Lambda}).
\end{equation*}
\end{prop}


Put for $\alpha \in \Sigma$
\begin{align*}
S_{\alpha,l}(\Lambda)
& = \{ \mu_1\cdots \mu_n \in S_l(\Lambda) \mid  
\mu_1 = \alpha \}, \\
S_{l+1}({}_\alpha \Lambda)
& = \{ \omega_1\cdots \omega_n \in S_{l+1}(\Lambda) \mid  
\alpha \omega_1\cdots \omega_n \in S_l(\Lambda) \}.
\end{align*}
For
$\omega=\omega_1\cdots\omega_n,
\gamma =\gamma_1\cdots \gamma_m
\in S_{l+1}({}_\alpha\Lambda)$,
we write
\begin{equation*}
\omega \underset{l,\alpha}{\sim}\gamma
\end{equation*}
if
$
\Gamma_l^-(\alpha\omega_1\cdots\omega_n)= \Gamma_l^-(\alpha\gamma_1\cdots \gamma_m).
$
We denote by 
$\Omega_\lambda^l({}_\alpha\Lambda)$
the equivalence classes of $S_{l+1}({}_\alpha\Lambda)$
under the equivalence relation
$\underset{l,\alpha}{\sim}$
with its discrete topology,
that is
$$
\Omega_\lambda^l({}_\alpha\Lambda)
= S_{l+1}({}_\alpha\Lambda) / \underset{l,\alpha}{\sim}.
$$
The equivalence class of
$\omega \in S_{l+1}({}_\alpha\Lambda)$
is denoted by
$[\omega]_{l,\alpha}$.
\begin{lem}
Assume that $\Lambda$ is $\lambda$-synchronizing.
For $\mu \in S_{\alpha,l}(\Lambda)$, 
there exists a word $\zeta =\zeta_1\cdots \zeta_n \in S_{l+1}({}_\alpha\Lambda)$
such that $\mu \underset{l}{\sim} \alpha \zeta$,
where 
$\alpha \zeta$ 
denotes 
$\alpha\zeta_1\cdots \zeta_n$.
\end{lem}
\begin{pf}
 For $\mu = \mu_1\cdots\mu_n \in S_{\alpha,l}(\Lambda)$ 
with
$\mu_1 =\alpha$, 
put
$k=n +l$.
As $\Lambda$ is $\lambda$-synchronizing,
there exists $\omega \in S_k(\Lambda)$
such that
$\mu \omega \in S_{k-n}(\Lambda)$.
Put
$\zeta =\mu_2 \cdots \mu_n \omega \in S_{k-(n-1)}(\Lambda) = S_{l+1}(\Lambda)$
so that
$\alpha \zeta = \mu \omega$
and
$\Gamma_l^-(\mu) = \Gamma_l^-(\mu \omega) =\Gamma_l^-(\alpha \zeta)$.
\end{pf}

Since
$
S_{l+2}({}_\alpha\Lambda)\subset S_{l+1}({}_\alpha\Lambda)
$
and
$\underset{l+1,\alpha}{\sim}$ 
implies
$\underset{l,\alpha}{\sim}$,
the map 
$
\iota_{l,l+1}^{\lambda(\Lambda)}:
\Omega_\lambda^{l+1}(\Lambda)
 \longrightarrow \Omega_\lambda^l(\Lambda)
$
induces a surjection
$$
\Omega_\lambda^{l+1}({}_\alpha\Lambda) \longrightarrow 
\Omega_\lambda^{l}({}_\alpha\Lambda)
\qquad \text{ for } l \in \Zp
$$
which gives rise to a compact Hausdorff space
$$
\underset{l}{\varprojlim} \Omega_\lambda^{l}({}_\alpha\Lambda)
$$
of a projective limit along the surjections.
We denote it by
 $\Omega_\lambda({}_\alpha\Lambda)$.

Let us consider the case for $\alpha=1$ in the above setting.
By the preceding lemma,
for $\mu \in S_{1,l}(\Lambda)$,
one may take $\zeta \in S_{l+1}({}_1\Lambda)$
such that
$\mu \underset{l}{\sim}1\zeta$.
One then sees the map
\begin{equation*}
[\mu]_l \in \Omega_\lambda^l(\Lambda_1) 
\longrightarrow
[\zeta]_{l,1} \in \Omega_\lambda^l({}_1\Lambda)
\end{equation*}
is well-defined.
We denote it by $s_l^\lambda$.
Conversely,
any element
$\omega_2\cdots\omega_n \in S_{l+1}({}_1\Lambda)$
yields 
$1\omega_2\cdots\omega_n \in S_{1,l}(\Lambda)$.
This correspondence yields the inverse of
$s^\lambda_l$,
so that 
$s^\lambda_l$ is bijective
and hence homeomorphic.
Since the sequence of the maps 
$\{ s^\lambda_l \}_{l \in \Zp}$
is  compatible to
$\iota$-maps
$$
\Omega_\lambda^{l+1}(\Lambda_1) \longrightarrow 
\Omega_\lambda^l(\Lambda_1)
\quad
\text{ and }
\quad
\Omega_\lambda^{l+1}({}_1\Lambda) \longrightarrow 
\Omega_\lambda^l({}_1\Lambda),
$$
they induce a homeomorphism
$$
\underset{l}{\varprojlim} s^\lambda_l :
\Omega_\lambda(\Lambda_1)
 \longrightarrow 
\Omega_\lambda({}_1\Lambda)
$$
denoted by
$s_\lambda.
$
By Proposition 5.8, we thus have
\begin{prop} 
The disjoint union
$
\Omega_\lambda({}_1\Lambda)
\sqcup
\Omega_\lambda(\Lambda)
$
is homeomorphic to
$\Omega_\lambda(\widetilde{\Lambda})
$
through the homeomorphism
$\varPhi_1 \circ s_\lambda^{-1} \sqcup \varPhi_d$.
\end{prop}
Denote by $\varPhi_\lambda$
the homeomorphism
$$
\varPhi_1 \circ s_\lambda^{-1} \sqcup \varPhi_d:
\Omega_\lambda({}_1\Lambda)
\sqcup
\Omega_\lambda(\Lambda)
\longrightarrow
\Omega_\lambda(\widetilde{\Lambda}).
$$
For $\alpha \in \Sigma$, denote by 
$\Omega_\lambda^l(\Lambda_\alpha)$
the $l$-past equivalence classes of 
$S_{\alpha,l}(\Lambda)$.
The map
$$
\mu_1 \cdots \mu_n \in S_{l+1}({}_\alpha\Lambda)
\longrightarrow 
\alpha \mu_1 \cdots \mu_n \in S_{\alpha,l}(\Lambda)
\subset S_l(\Lambda)
$$
induces a map
\begin{equation*}
\lambda_l(\alpha):
[\mu ]_{l,\alpha} \in \Omega_\lambda^l({}_\alpha\Lambda) 
\longrightarrow
[\alpha \mu ]_l \in \Omega_\lambda^l(\Lambda_\alpha)
\subset
\Omega_\lambda^l(\Lambda). 
\end{equation*}
The family 
$\{ \lambda_l(\alpha) \}_{l \in \Zp}$
yields a continuous map
:
$
\Omega_\lambda({}_\alpha\Lambda) 
\longrightarrow
\Omega_\lambda(\Lambda) 
$
which we write
$
\lambda_\Lambda(\alpha).
$
Let
${\Bbb Z}_{\lambda(\Lambda)}$
and
${\Bbb Z}_{\lambda({}_\alpha\Lambda)}$
be the abelian groups 
$C(\Omega_\lambda(\Lambda),{\Bbb Z})$
and
$C(\Omega_\lambda({}_\alpha\Lambda),{\Bbb Z})$
of all ${\Bbb Z}$-valued continuous functions on the compact Hausdorff spaces
$\Omega_\lambda(\Lambda)$
and
$\Omega_\lambda({}_\alpha\Lambda)$
respectively.
The map
$
\lambda_\Lambda(\alpha)
$
induces a homomorphism of abelian groups
$$
\lambda_\Lambda(\alpha)^*:
{\Bbb Z}_\lambda(\Lambda)
\longrightarrow
{\Bbb Z}_\lambda({}_\alpha\Lambda)
$$
such that 
$$
\lambda_\Lambda(\alpha)^*(f)(t) = f(\lambda_\Lambda(\alpha)(t))
\quad
\text{ for }
\qquad f \in {\Bbb Z}_\lambda(\Lambda), 
t \in \Omega_\lambda({}_\alpha\Lambda).
$$
Let
${}_\alpha\Omega_\lambda^{l+1}(\Lambda)$
be 
the $l+1$-past equivalence classes of $S_{l+1}({}_\alpha\Lambda)$
so that
$$
{}_\alpha\Omega_\lambda^{l+1}(\Lambda) 
= \{ [\mu]_{l+1}\in \Omega_\lambda^{l+1}(\Lambda) \mid
\alpha \in \Gamma^-_1(\mu)\} 
\subset \Omega_\lambda^{l+1}(\Lambda).
$$
The surjective map
$$
\iota_{l,l+1}^{\lambda(\Lambda)}:
\Omega_\lambda^{l+1}(\Lambda) \longrightarrow
\Omega_\lambda^l(\Lambda)
$$
works for the restriction
$$
{}_\alpha\Omega_\lambda^{l+1}(\Lambda) \longrightarrow
{}_\alpha\Omega_\lambda^l(\Lambda)
$$
which yields a compact Hausdorff space 
$$
{}_\alpha\Omega_\lambda(\Lambda)
=\underset{l}{\varprojlim} \ {}_\alpha\Omega_\lambda^l(\Lambda)
$$
by taking a projective limit along the above restrictions.
We may regard 
$
{}_\alpha\Omega_\lambda(\Lambda)
$
as a clopen subset of
$
\Omega_\lambda(\Lambda).
$

For 
$\mu,\nu \in S_{l+1}({}_\alpha\Lambda)$,
the condition
$[\mu]_{l+1} = [\nu]_{l+1}$
implies 
$
[\mu]_{l,\alpha} = 
[\nu]_{l,\alpha}
$
so that the map
$$
e_{(\alpha)}^l: 
[\mu]_{l+1} \in {}_\alpha\Omega_\lambda^{l+1}(\Lambda)
\longrightarrow
[\mu]_{l,\alpha} \in \Omega_\lambda^{l}({}_\alpha\Lambda)
\qquad
\text{ for }
\quad
\mu \in S_{l+1}({}_\alpha\Lambda)
$$
is well-defined and surjective.
Since the maps
$
e_{(\alpha)}^l: 
{}_\alpha\Omega_\lambda^{l+1}(\Lambda)
\longrightarrow
 \Omega_\lambda^{l}({}_\alpha\Lambda),
 l \in \Zp
$
are compatible to the surjections
$$
\iota_{l,l+1}^{\lambda(\Lambda)}|_{ {}_\alpha\Omega_\lambda^{l+1}(\Lambda) }:
{}_\alpha\Omega_\lambda^{l+1}(\Lambda) 
\longrightarrow
{}_\alpha\Omega_\lambda^l(\Lambda),
\qquad
\iota_{l,l+1}^{\lambda(\Lambda)} :
\Omega_\lambda^{l+1}({}_\alpha\Lambda) 
\longrightarrow
\Omega_\lambda^l({}_\alpha\Lambda),
$$
we have a continuous surjection
$$
e_{(\alpha)}^{\lambda(\Lambda)} 
= \varprojlim e_{(\alpha)}^l :{}_\alpha\Omega_\lambda(\Lambda)
\longrightarrow
\Omega_\lambda({}_\alpha\Lambda).
$$
We then define a homomorphism of abelian groups
$$
\epsilon_{(\alpha)}^{\lambda(\Lambda)} :
{\Bbb Z}_{\lambda({}_\alpha\Lambda)}
\longrightarrow
{\Bbb Z}_{\lambda(\Lambda)}
$$
by setting
\begin{equation*}
\epsilon_{(\alpha)}^{\lambda(\Lambda)}(f)(x) =
\begin{cases}
f(e_{(\alpha)}^{\lambda(\Lambda)}(x)) 
  & \text{ if } x \in {}_\alpha\Omega_\lambda(\Lambda),\\
0 & \text{ otherwise }
\end{cases}
\end{equation*}
 for
$
f \in C(\Omega_\lambda({}_\alpha\Lambda), {\Bbb Z}),
x \in \Omega_\lambda(\Lambda).
$
Thus we have a family of endomorphisms of abelian groups
$$
\epsilon_{(\alpha)}^{\lambda(\Lambda)} \circ \lambda_\Lambda(\alpha)^* :
{\Bbb Z}_{\lambda(\Lambda)} 
\longrightarrow
{\Bbb Z}_{\lambda(\Lambda)}.
$$
Let us define an endomorphism 
$\lambda(\Lambda)$
of  
${\Bbb Z}_{\lambda(\Lambda)}$
by setting:
$$
\lambda(\Lambda) =
\sum_{\alpha=1}^{N} 
\epsilon_{(\alpha)}^{\lambda(\Lambda)} \circ \lambda_\Lambda(\alpha)^* :
{\Bbb Z}_{\lambda(\Lambda)} 
\longrightarrow
{\Bbb Z}_{\lambda(\Lambda)}.
$$ 
By definition of the $\lambda$-synchronizing K-groups,
we see
\begin{prop}
\hspace{3cm}
\begin{enumerate}
\renewcommand{\labelenumi}{(\roman{enumi})}
\item
$
K_0^\lambda(\Lambda) 
= {\Bbb Z}_{\lambda(\Lambda)}
 /(\id -\lambda(\Lambda)){\Bbb Z}_{\lambda(\Lambda)}.
$
\item
$
K_1^\lambda(\Lambda) = {\Ker }(\id -\lambda(\Lambda))
$ 
  in  
$
{\Bbb Z}_{\lambda(\Lambda)}.
$
\end{enumerate}
\end{prop}
\begin{pf}
Let
$\{ v_i^l \}_{i=1}^{m_\lambda(l)}$
be the vertex set $V^{\lambda(\Lambda)}_l$ of the 
$\lambda$-synchronizing $\lambda$-graph system
$\LLL$.
As we identify
$\Omega_\lambda^l(\Lambda)$ with $V^{\lambda(\Lambda)}_l$,
we may write
$v_i^l = [\mu^l(i)]_l$ for some 
$\mu^l(i) \in S_l(\Lambda), i=1,\dots,m_\lambda(l).$ 
For 
$
v_i^l \in V^{\lambda(\Lambda)}_l, 
v_j^{l+1} \in  V^{\lambda(\Lambda)}_{l+1}
$
and
$\alpha \in \Sigma$,
we set
$A^{\lambda(\Lambda)}_{l,l+1}(i,\alpha,j) =1$
if there exists an labeled edge labeled $\alpha$
which starts at $v_i^l$ 
and ends at $v_j^{l+1}$,
and 
$A^{\lambda(\Lambda)}_{l,l+1}(i,\alpha,j) =0$,
otherise.
Hence we have
$A^{\lambda(\Lambda)}_{l,l+1}(i,j)
= \sum_{\alpha=1}^N A^{\lambda(\Lambda)}_{l,l+1}(i,\alpha,j).
$
We then have for $j=1,\dots,m_\lambda(l+1)$
\begin{align*}
& \{ [\alpha\mu^{l+1}(j)]_l \in V_l^{\lambda(\Lambda)} \mid   
  \alpha\in \Gamma_1^-(\mu^{l+1}(j)) \}\\
=
& \{ [\mu^l(i)]_l \in V_l^{\lambda(\Lambda)} \mid  
A^{\lambda(\Lambda)}_{l,l+1}(i,\alpha,j) =1, i=1,\dots,m_\lambda(l), \alpha \in \Sigma \}.
\end{align*} 
For $ f \in C(\Omega_\lambda^l(\Lambda),{\Bbb Z})$
and $v_j^{l+1} \in  V^{\lambda(\Lambda)}_{l+1}, 
j=1,\dots,m_\lambda(l+1)$
it follows that
\begin{align*}
\lambda(\Lambda)(f) (v_j^{l+1})
&= \sum_{\alpha=1}^{N} 
\epsilon_{(\alpha)}^{\lambda(\Lambda)} \circ \lambda_\Lambda(\alpha)^*(f)
([\mu^{l+1}(j)]_{l+1})\\
&= \sum_{\alpha\in \Gamma_1^-(\mu^{l+1}(j))}f([\alpha\mu^{l+1}(j)]_l)\\
&= \sum_{i=1}^{m_\lambda(l)}\sum_{\alpha=1}^N 
A^{\lambda(\Lambda)}_{l,l+1}(i,\alpha,j) f([\mu^l(i)]_l)\\
&= \sum_{i=1}^{m_\lambda(l)}
A^{\lambda(\Lambda)}_{l,l+1}(i,j) f([\mu^l(i)]_l)\\
&= \sum_{i=1}^{m_\lambda(l)}
f(v_i^l){}^t\negthinspace A^{\lambda(\Lambda)}_{l,l+1}(j,i).
\end{align*}
The group
$C(\Omega_\lambda^l(\Lambda),{\Bbb Z})$
is identified with
${\Bbb Z}^{m_\lambda(l)}$
through the map
$$
f \in C(\Omega_\lambda^l(\Lambda),{\Bbb Z})
\longrightarrow
[f(v_i^l)]_{i=1}^{m_\lambda(l)}\in {\Bbb Z}^{m_\lambda(l)}.
$$
Hence the homomorphism
$$
\lambda(\Lambda) |_{C(\Omega_\lambda^l(\Lambda),{\Bbb Z})}:
C(\Omega_\lambda^l(\Lambda),{\Bbb Z})
\longrightarrow
C(\Omega_\lambda^{l+1}(\Lambda),{\Bbb Z})
$$
is identified with the matrix
$$
{}^tA^{\lambda(\Lambda)}_{l,l+1}:
{\Bbb Z}^{m_\lambda(l)} 
\longrightarrow 
{\Bbb Z}^{m_\lambda(l+1)}. 
$$
Therefore by definition 
of $K_i^\lambda(\Lambda), i=0,1$,
we have  desired formulae.
\end{pf}

\begin{lem}
There exists a natural identification
$$
\Omega_\lambda^l(\widetilde{\Lambda}_1) 
={}_0\Omega_\lambda^l(\widetilde{\Lambda}) 
\text{ for all }
l \in \Zp
\quad
\text{ so that }
\quad
\Omega_\lambda(\widetilde{\Lambda}_1) 
={}_0\Omega_\lambda(\widetilde{\Lambda}). 
$$
\end{lem}
\begin{pf}
Since for 
$\mu = \mu_1\cdots\mu_n \in S_l(\widetilde{\Lambda})$
the condition $\mu_1 =1$
is equivalent to the condition
$
0\mu_1\cdots\mu_n \in S_{0,l-1}(\widetilde{\Lambda}),
$
one has
$
S_{1,l}(\widetilde{\Lambda}) 
=S_l({}_0\widetilde{\Lambda}) 
$
so that
$
\Omega_\lambda^l(\widetilde{\Lambda}_1) 
={}_0\Omega_\lambda^l(\widetilde{\Lambda}).
$ 
\end{pf}
\begin{lem}
\hspace{3cm}
\begin{enumerate}
\renewcommand{\labelenumi}{(\roman{enumi})}
\item
$
\varPhi_d^{-1}
\circ \lambda_{\widetilde{\Lambda}}(0)
\circ e_{(0)}^{\lambda(\widetilde{\Lambda})}
\circ \varPhi_1
\circ s_\lambda^{-1} 
= \lambda_\Lambda(1)
$
on 
$\Omega_\lambda({}_1\Lambda)$.
\item
$
\varPhi_\lambda^{-1}
\circ \lambda_{\widetilde{\Lambda}}(\alpha)
\circ e_{(\alpha)}^{\lambda(\widetilde{\Lambda})}
\circ \varPhi_d
=
{\begin{cases}
e_{(1)}^{\lambda(\Lambda)}
\text{ on }
{}_1\Omega_\lambda(\Lambda)
& \text{ for } \alpha =1,\\
\lambda_\Lambda(\alpha) 
\circ e_{(\alpha)}^{\lambda(\Lambda)}
\text{ on }
{}_\alpha\Omega_\lambda(\Lambda)
& \text{ for } \alpha=2,3,\cdots N.
\end{cases}}
$
\end{enumerate}
\end{lem}
\begin{pf}
(i)
We will prove that the equality
\begin{equation*}
\varPhi_d \circ \lambda_\Lambda(1) \circ s_\lambda
=
\lambda_{\widetilde{\Lambda}}(0)
\circ e_{(0)}^{\lambda(\widetilde{\Lambda})}
\circ \varPhi_1
\quad
\text{ on }
\Omega_\lambda(\Lambda_1)
\end{equation*}
holds.
For 
$\mu \in S_{1,l}(\Lambda)$,
take
$\zeta \in S_{l+1}({}_1\Lambda)$
such that
$\mu \underset{l}{\sim} 1\zeta$.
It then follows that
\begin{equation*}
(\xi_l^\lambda \circ \lambda_l(1) \circ s_l^\lambda)([\mu]_l)
 = ( \xi_l^\lambda \circ \lambda_l(1))([\zeta]_{l,1}) 
 =  \xi_l^\lambda ([\mu]_l) 
 = [ \xi^B(\mu) ]_l.
\end{equation*}
On the other hand,
\begin{align*}
  & (\lambda_l(0)\circ e_{(0)}^l\circ \varphi_{l+1}^\lambda )([\mu]_l) \\
= & (\lambda_l(0)\circ e_{(0)}^l)([\varphi^B(\mu)]_{l+1})
= \lambda_l(0)([\varphi^B(\mu)]_{l,0})
=  [0\varphi^B(\mu)]_l.
\end{align*}
As 
$\varphi^B(\mu) = 1\xi^B(\mu_2\cdots\mu_n)$,
one has
$$
0\varphi^B(\mu)= 0 1\xi^B(\mu_2\cdots\mu_n) = \xi^B(\mu)
$$
so that we have
$
\varPhi_d \circ \lambda_\Lambda(1) \circ s_\lambda 
=\lambda_{\widetilde{\Lambda}}(0)
\circ e_{(0)}^{\lambda(\widetilde{\Lambda})}
\circ \varPhi_1.
$

(ii)
For $\alpha=1$, 
the map
$\lambda_{\widetilde{\Lambda}}(1)$
is defined from 
$\Omega_\lambda({}_1\widetilde{\Lambda})$
to
$\Omega_\lambda(\widetilde{\Lambda}_1)$.
As
$
\varPhi_{\lambda} |_{\Omega_\lambda({}_1\Lambda)}
=\varPhi_1 \circ s_\lambda^{-1}:
\Omega_\lambda({}_1\Lambda)
\longrightarrow 
\Omega_\lambda(\widetilde{\Lambda}_1),
$
one sees 
$\varPhi_\lambda^{-1} |_{\Omega_\lambda(\widetilde{\Lambda}_1)}
=s_\lambda \circ \varPhi_1^{-1}.
$
For 
$[\gamma]_l \in {}_1\Omega_\lambda^l(\Lambda)
$
with
$\gamma \in S_l({}_1\Lambda)$,
it follows that
\begin{equation*}
(\lambda_{l-1}(1)
\circ e_{(1)}^{l-1}
\circ \xi_l^\lambda) ([\gamma]_l)
 =(\lambda_{l-1}(1)\circ e_{(1)}^{l-1}) ([\tilde{\gamma}]_l) 
 =\lambda_{l-1}(1) ([\tilde{\gamma}]_{l-1,1}) 
 =[1\tilde{\gamma}]_{l-1}.
\end{equation*}
Since
\begin{equation*}
(\varphi_l^\lambda \circ (s_{l-1}^\lambda )^{-1}\circ  e_{(1)}^{l-1}) ([\gamma]_l)
 = (\varphi_l^\lambda \circ (s_{l-1}^\lambda)^{-1})([ \gamma]_{l-1,1})
 = \varphi_l^\lambda ([1 \gamma]_{l-1})
 =[1\tilde{\gamma}]_l.
\end{equation*}
As
$\iota_{l-1,l}^{\lambda(\widetilde{\Lambda})}([1\tilde{\gamma}]_l) 
=[1\tilde{\gamma}]_{l-1},
$
we have 
$
 \lambda_{\widetilde{\Lambda}}(1)
\circ e_{(1)}^{\lambda(\widetilde{\Lambda})}
\circ \varPhi_d 
=
\varPhi_\lambda \circ
 e_{(1)}^{\lambda(\Lambda)}.
$

For $\alpha = 2,3,\cdots, N$, 
the map
$\lambda_{\widetilde{\Lambda}}(\alpha)$
is defined from 
$\Omega_\lambda({}_\alpha\widetilde{\Lambda})$
to
$\Omega_\lambda(\widetilde{\Lambda}_\alpha)$.
As
$
\Omega_\lambda(\widetilde{\Lambda}_\alpha)
$
is contained in 
$
\Omega_\lambda(\widetilde{\Lambda})
\backslash
\Omega_\lambda(\widetilde{\Lambda}_1),
$
one has 
$
\varPhi_\lambda^{-1}\circ \lambda_{\widetilde{\Lambda}}(\alpha)
= \varPhi_d^{-1}\circ \lambda_{\widetilde{\Lambda}}(\alpha).
$
For 
$[\gamma]_l \in {}_\alpha\Omega^l_\lambda(\Lambda)
$
with
$\gamma \in S_l({}_\alpha\Lambda)$,
it follows that
\begin{equation*}
(\lambda_{l-1}(\alpha)\circ e_{(\alpha)}^{l-1}\circ \xi_l^\lambda) ([\gamma]_l)
 =(\lambda_{l-1}(\alpha)\circ e_{(\alpha)}^{l-1}) ([\tilde{\gamma}]_l) 
 =\lambda_{l-1}(\alpha) ([\tilde{\gamma}]_{l-1,\alpha}) 
 = [\alpha\tilde{\gamma}]_{l-1}.
\end{equation*}
On the other hand,
\begin{equation*}
(\xi_{l-1}^\lambda\circ \lambda_{l-1}(\alpha) 
\circ e_{(\alpha)}^{l-1}) ([\gamma]_l)
=(\xi_{l-1}^\lambda\circ \lambda_{l-1}(\alpha) ) ([\gamma]_{l-1,\alpha})
=\xi_{l-1}^\lambda([\alpha \gamma]_{l-1})
= [\alpha\tilde{\gamma}]_{l-1}
\end{equation*}
so that we have
$
\lambda_{\widetilde{\Lambda}}(\alpha)
\circ e_{(\alpha)}^{\lambda(\widetilde{\Lambda})}
\circ \varPhi_d 
=
\varPhi_\lambda \circ
\lambda_\Lambda(\alpha) 
\circ e_{(\alpha)}^{\lambda(\Lambda)}. 
$
\end{pf}
Consider an endomorphism on the group
$
{\Bbb Z}_{\lambda({}_1\Lambda)}
\oplus {\Bbb Z}_{\lambda(\Lambda)}
$
defined by
$$
(g,h) \in {\Bbb Z}_{\lambda({}_1\Lambda)}
\oplus {\Bbb Z}_{\lambda(\Lambda)}
\longrightarrow
(\lambda_\Lambda(1)^*(h),
\epsilon_{(1)}^{\lambda(\Lambda)}(g) 
+ \sum_{\alpha=2}^{N} (\epsilon_{(\alpha)}^{\lambda(\Lambda)} \circ \lambda_\Lambda(\alpha)^*) (h))
 \in {\Bbb Z}_{\lambda({}_1\Lambda)}\oplus {\Bbb Z}_{\lambda(\Lambda)}
$$
that is represented  by
$$
A_{\lambda(\Lambda)} =
\begin{bmatrix}
0,  & \lambda_\Lambda(1)^* \\
\epsilon_{(1)}^{\lambda(\Lambda)}, 
& \sum_{\alpha=2}^{N} \epsilon_{(\alpha)}^{\lambda(\Lambda)} \circ \lambda_\Lambda(\alpha)^*
\end{bmatrix}
\qquad
\text{on }
\quad
{\Bbb Z}_{\lambda({}_1\Lambda)}\oplus {\Bbb Z}_{\lambda(\Lambda)}.
$$
For the subshift $\widetilde{\Lambda}$,
we will similarly formulate the
endomorphism 
$
\lambda({\tilde{\Lambda}})
$ 
on
${\Bbb Z}_{\lambda(\tilde{\Lambda})}$
as in the following way.
Recall that for $\alpha\in \widetilde{\Sigma}$ 
$$
S_{l+1}({}_\alpha\widetilde{\Lambda})
 = \{ \gamma_1\cdots\gamma_m 
 \in S_{l+1}(\widetilde{\Lambda}) \mid
 \alpha\gamma_1\cdots\gamma_m \in  S_l(\widetilde{\Lambda})
 \}.
$$
The sets 
${}_\alpha\Omega_\lambda^{l+1}(\widetilde{\Lambda})$
and
$\Omega_\lambda^l({}_\alpha\widetilde{\Lambda})$
are defined by the equivalence classes
$
S_{l+1}({}_\alpha\widetilde{\Lambda}) / \underset{l+1}{\sim}
$
and
$
S_{l+1}({}_\alpha\widetilde{\Lambda}) / \underset{l,\alpha}{\sim}
$
respectively.
The compact Hausdorff spaces 
${}_\alpha\Omega_\lambda(\widetilde{\Lambda})$
and
$\Omega_\lambda({}_\alpha\widetilde{\Lambda})$
are defined by the projective limits
$\underset{l}{\varprojlim} {}_\alpha\Omega_\lambda^{l+1}(\widetilde{\Lambda})$
and
$\underset{l}{\varprojlim} \Omega_\lambda^l({}_\alpha\widetilde{\Lambda})$
respectively.
The maps
$\tilde{e}_{(\alpha)}^l:{}_\alpha\Omega_\lambda^{l+1}(\widetilde{\Lambda})
\longrightarrow
\Omega_\lambda^l({}_\alpha\widetilde{\Lambda})
$
and
$\tilde{\lambda}_{l}(\alpha):\Omega_\lambda^{l}({}_\alpha\widetilde{\Lambda})
\longrightarrow
\Omega_\lambda^l(\widetilde{\Lambda})
$
are defined by 
\begin{align*}
\tilde{e}_{(\alpha)}^l([\gamma]_l) & = [\gamma]_{l,\alpha}, 
\qquad \text{ for } \gamma \in S_{l+1}({}_\alpha\widetilde{\Lambda}),\\
\tilde{\lambda}_{l}(\alpha)([\mu]_{l,\alpha}) & = [\alpha\mu]_l, 
\qquad \text{ for } \mu \in S_{l+1}({}_\alpha\widetilde{\Lambda}).
\end{align*}
Set
\begin{align*}
e_{(\alpha)}^{\lambda(\widetilde{\Lambda})}
& = \underset{l}{\varprojlim} \tilde{e}_{(\alpha)}^l:
{}_\alpha\Omega_\lambda(\widetilde{\Lambda})
\longrightarrow
\Omega_\lambda({}_\alpha\widetilde{\Lambda}),\\
\lambda_{\widetilde{\Lambda}}(\alpha)
& = \underset{l}{\varprojlim}\tilde{\lambda}_{l}(\alpha):
\Omega_\lambda({}_\alpha\widetilde{\Lambda})
\longrightarrow
\Omega_\lambda(\widetilde{\Lambda}).
\end{align*}
The map
$\epsilon_{(\alpha)}^{\lambda(\widetilde{\Lambda})}:
{\Bbb Z}_{\lambda({}_\alpha\widetilde{\Lambda})} \longrightarrow
{\Bbb Z}_{\lambda(\widetilde{\Lambda})}
$
is defined 
by 
\begin{equation*}
\epsilon_{(\alpha)}^{\lambda(\widetilde{\Lambda})}(f)(x)
= 
\begin{cases}
f(e_{(\alpha)}^{\lambda(\widetilde{\Lambda})}(x)) & 
\text{ if } x \in {}_\alpha\Omega_\lambda(\widetilde{\Lambda}),\\
0 & 
\text{ otherwise, }
\end{cases}
\end{equation*}
for 
$
f \in {\Bbb Z}_{\lambda({}_\alpha\widetilde{\Lambda})},
x \in \Omega_\lambda(\widetilde{\Lambda}).
$
The endomorphism
$
\lambda({\tilde{\Lambda}})
$
on
${\Bbb Z}_{\lambda(\tilde{\Lambda})}$
is defined by 
\begin{equation*}
\lambda({\tilde{\Lambda}})
= \sum_{\alpha=0}^N \epsilon_{(\alpha)}^{\lambda(\widetilde{\Lambda})} \circ 
\lambda_{\widetilde{\Lambda}}(\alpha)^*.
\end{equation*}
The homeomorphism 
$
\varPhi_\lambda: 
\Omega_\lambda({}_1\Lambda)\sqcup \Omega_\lambda(\Lambda)
\longrightarrow \Omega_\lambda(\widetilde{\Lambda}) 
$ 
naturally
yields an isomorphism $\varPhi_\lambda^*$ 
of abelian groups
from
${\Bbb Z}_{\lambda(\widetilde{\Lambda})}$
to
${\Bbb Z}_{\lambda({}_1\Lambda)}
\oplus 
{\Bbb Z}_{\lambda(\Lambda)}.$
We then have
\begin{lem}
 $\lambda(\tilde{\Lambda})
  = \varPhi_\lambda^{*-1} 
  \circ A_{\lambda(\Lambda)} \circ \varPhi_\lambda^*$
on 
${\Bbb Z}_{\lambda(\widetilde{\Lambda})}.$
\end{lem}
\begin{pf}
The proof is completely similar to the proof of
\cite[Lemma 2.8]{2001ETDS}.
We will give the proof for the sake of completeness.
For
$
(g,h) \in 
{\Bbb Z}_{\lambda({}_1\Lambda)}\oplus {\Bbb Z}_{\lambda(\Lambda)}
$  
and 
$
\gamma \in \Omega_\lambda({}_1\Lambda)\sqcup \Omega_\lambda(\Lambda),
$
we see
$$
[(\varPhi_\lambda^* \circ \lambda(\tilde{\Lambda})
\circ
\varPhi_\lambda^{*-1}) (g,h)](\gamma)
=(\lambda(\tilde{\Lambda})(
\varPhi_\lambda^{*-1} (g,h))(\varPhi_\lambda(\gamma)).
$$
We have two cases.

\noindent
Case 1: $\gamma \in \Omega_\lambda({}_1\Lambda)$.

Since 
$
\varPhi_1\circ s_\lambda^{-1}: 
\Omega_\lambda({}_1\Lambda) \longrightarrow 
\Omega_\lambda(\widetilde{\Lambda}_1)
$
and
$
\Omega_\lambda(\widetilde{\Lambda}_1) 
= {}_0\Omega_\lambda(\widetilde{\Lambda}),
$
one has
$\varPhi_\lambda(\gamma) 
= \varPhi_1 \circ s_\lambda^{-1}(\gamma) \in {}_0\Omega_\lambda(\widetilde{\Lambda}).
$
We have 
$
\epsilon_{(\alpha)}^{\lambda(\widetilde{\Lambda})} 
\circ \lambda_{\widetilde{\Lambda}}(\alpha)^* 
\circ \varPhi_\lambda^{*-1} (g,h) (\varPhi_\lambda (\gamma) ) =  0
$
for $\alpha \ne 0$. 
By the preceding lemma, one sees 
\begin{align*}
[(\lambda(\widetilde{\Lambda}) \circ \varPhi_\lambda^{*-1}) (g,h)] 
(\varPhi_\lambda(\gamma)) 
= & \sum_{\alpha =0}^{N}
\epsilon_{(\alpha)}^{\lambda(\widetilde{\Lambda})} \circ \lambda_{\widetilde{\Lambda}}(\alpha)^* 
    \circ \varPhi_\lambda^{*-1} (g,h) (\varPhi_\lambda (\gamma) ) \\
= &
\epsilon_{(0)}^{\lambda(\widetilde{\Lambda})} 
\circ \lambda_{\widetilde{\Lambda}}(0)^* 
\circ \varPhi_\lambda^{*-1} (g,h) (\varPhi_1 \circ s_\lambda^{-1} (\gamma) ) \\
= &
(g,h)(\varPhi_\lambda^{-1}\circ \lambda_{\widetilde{\Lambda}}(0)\circ 
e_{(0)}^{\lambda(\widetilde{\Lambda})}
\circ \varPhi_1\circ s_\lambda^{-1} )(\gamma)  \\
= & (g,h)(\lambda_\Lambda(1)(\gamma))\\
= & h(\lambda_\Lambda(1)(\gamma)) = (\lambda_\Lambda(1)^*(h))(\gamma)).
\end{align*}

Case 2: $\gamma \in \Omega_\lambda(\Lambda)$.

As 
$
\varPhi_d (\gamma) \not\in\Omega_\lambda(\widetilde{\Lambda}_1)
= {}_0\Omega_\lambda(\widetilde{\Lambda}),
$
one sees that
\begin{equation*}
\epsilon_{(0)}^{\lambda(\widetilde{\Lambda})} 
\circ \lambda_{\widetilde{\Lambda}}(0)^* 
\circ \varPhi_\lambda^{*-1} (g,h) (\varPhi_d (\gamma) )=0.
\end{equation*}
It follows that
\begin{equation*}
[\lambda({\tilde{\Lambda}}) \circ \varPhi_\lambda^{*-1} (g,h)] (\varPhi(\gamma)) 
=  \sum_{\alpha=1}^{N}
\epsilon_{(\alpha)}^{\lambda(\widetilde{\Lambda})} 
\circ \lambda_{\widetilde{\Lambda}}(\alpha)^* 
\circ \varPhi_\lambda^{*-1} (g,h) (\varPhi_d (\gamma) ). 
\end{equation*}
It is easy to see that
$
\varPhi_d (\gamma)
\in {}_\alpha\Omega_\lambda(\widetilde{\Lambda})
$
if and only if
$
\gamma
\in {}_\alpha\Omega_\lambda(\Lambda)
$
for 
$\alpha = 1,2,\dots,N$.
Hence 
$
\epsilon_{(\alpha)}^{\lambda(\widetilde{\Lambda})} 
\circ \lambda_{\widetilde{\Lambda}}(\alpha)^* 
    \circ \varPhi_\lambda^{*-1} (g,h) (\varPhi_d (\gamma) ) =0
$
for
$\gamma \not\in {}_\alpha\Omega_\lambda(\Lambda).
$
If 
$
\gamma \in {}_1\Omega_\lambda(\Lambda),
$
 we have by the preceding lemma
\begin{align*}
\epsilon_{(1)}^{\lambda(\widetilde{\Lambda})} 
\circ \lambda_{\widetilde{\Lambda}}(1)^* 
\circ \varPhi_\lambda^{*-1} (g,h) (\varPhi_d (\gamma) ) 
= & (g,h)(\varPhi_\lambda^{-1}\circ \lambda_{\widetilde{\Lambda}}(1)
\circ e_{(1)}^{\lambda(\widetilde{\Lambda})}\circ \varPhi_d (\gamma))\\
= & (g,h)(e_{(1)}^{\lambda(\Lambda)}(\gamma)) \\
= & g(e_{(1)}^{\lambda(\Lambda)}(\gamma)) = \epsilon_{(1)}^{\lambda(\Lambda)}(g)(\gamma)).
\end{align*}
If 
$
\gamma \in {}_\alpha\Omega_\lambda(\Lambda)
$
for 
$\alpha=2,3,\dots, N$, 
we have  by the preceding lemma,
\begin{align*}
\epsilon_{(\alpha)}^{\lambda(\widetilde{\Lambda})} 
\circ \lambda_{\widetilde{\Lambda}}(\alpha)^* 
\circ \varPhi_\lambda^{*-1} (g,h) (\varPhi_d (\gamma) ) 
= & (g,h)(\varPhi_\lambda^{-1} 
\circ \lambda_{\widetilde{\Lambda}}(\alpha)
\circ e_{(\alpha)}^{\lambda(\widetilde{\Lambda})}\circ \varPhi_d (\gamma) )\\
= & h(\varPhi_\lambda^{-1}
\circ \lambda_{\widetilde{\Lambda}}(\alpha)
\circ e_{(\alpha)}^{\lambda(\widetilde{\Lambda})}\circ \varPhi_d (\gamma)) \\
= & h( \lambda_{\Lambda}(\alpha)
\circ e_{(\alpha)}^{\lambda(\Lambda)} (\gamma) )\\
= &[( \epsilon_{(\alpha)}^{\lambda(\Lambda)} \circ \lambda_\Lambda(\alpha)^*)(h)] (\gamma).
\end{align*}
Therefore  we conclude
\begin{equation*}
[\lambda({\tilde{\Lambda}}) 
\circ \varPhi_\lambda^{*-1} (g,h)] (\varPhi_\lambda(\gamma)) 
=  
\epsilon_{(1)}^{\lambda(\Lambda)}(g)(\gamma)
+
\sum_{\alpha=2}^N[( \epsilon_{(\alpha)}^{\lambda(\Lambda)} 
\circ \lambda_\Lambda(\alpha)^*)(h)] (\gamma)
\end{equation*}
so that
\begin{equation*}
\varPhi_\lambda^{*} 
\circ\lambda(\tilde{\Lambda})  
\circ \varPhi_\lambda^{*-1} 
  = A_{\lambda(\Lambda)}. 
\end{equation*}
\end{pf}

Our goal is established by showing the following isomorphisms of groups:
\begin{lem}
\hspace{3cm}
\begin{enumerate}
\renewcommand{\labelenumi}{(\roman{enumi})}
\item
$
{\Bbb Z}_{\lambda(\Lambda)}
 /(\id -\lambda(\Lambda)){\Bbb Z}_{\lambda(\Lambda)}
 \cong
 ({\Bbb Z}_{\lambda({}_1\Lambda)} 
\oplus {\Bbb Z}_{\lambda(\Lambda)} ) /
  (\id - A_{\lambda(\Lambda)}) 
  ({\Bbb Z}_{\lambda({}_1\Lambda)} \oplus{\Bbb Z}_{\lambda(\Lambda)})
$
\item 
$
{\Ker }(\id -\lambda(\Lambda))
$ 
in 
$
{\Bbb Z}_{\lambda(\Lambda)}
\cong
{\Ker }(\id - A_{\lambda(\Lambda)})
$
 in 
$
{\Bbb Z}_{\lambda({}_1\Lambda)} \oplus{\Bbb Z}_{\lambda(\Lambda)}.
$
\end{enumerate}
\end{lem}
The proof given here is similar to the proof of \cite[Lemma 2.9]{2001ETDS}
that is basically due to the original proof of Bowen-Franks in \cite{BF}.
\begin{pf}
(i)
Define a homomorphism
$$
\delta: (g,h) \in {\Bbb Z}_{\lambda({}_1\Lambda)} \oplus{\Bbb Z}_{\lambda(\Lambda)}
\longrightarrow
[\epsilon_{(1)}^{\lambda(\Lambda)}(g) +h] \in  
{\Bbb Z}_{\lambda(\Lambda)} /(\id -\lambda_{(\Lambda)}){\Bbb Z}_{\lambda(\Lambda)}.
$$
It is clear that $\delta$ is surjective.
We will show that
$$
\Ker(\delta) 
= (\id - A_{\lambda(\Lambda)}) 
({\Bbb Z}_{\lambda({}_1\Lambda)} \oplus{\Bbb Z}_{\lambda(\Lambda)}).
$$
For 
$ 
(g,h) \in {\Bbb Z}_{\lambda({}_1\Lambda)} \oplus{\Bbb Z}_{\lambda(\Lambda)},
$
it follows that
\begin{align*}
\delta((\id - A_{\lambda(\Lambda)}))(g,h)
= & \delta ((g - \lambda_\Lambda(1)^*(h), 
    h - \epsilon_{(1)}^{\lambda(\Lambda)}(g) 
      - \sum_{\alpha=2}^{N} \epsilon_{(\alpha)}^{\lambda(\Lambda)}\circ\lambda_\Lambda(\alpha)^*(h))\\
= & [\epsilon_{(1)}^{\lambda(\Lambda)}  ( g - \lambda_\Lambda(1)^*(h) ) 
  + h - \epsilon_{(1)}^{\lambda(\Lambda)}(g) 
  - \sum_{\alpha=2}^{N} \epsilon_{(\alpha)}^{\lambda(\Lambda)} \circ\lambda_\Lambda(\alpha)^*(h)]\\
= & [ h - \sum_{\alpha=1}^{N} \epsilon_{(\alpha)}^{\lambda(\Lambda)} \circ\lambda_\Lambda(\alpha)^*(h)]\\
= & [(\id - \lambda(\Lambda))(h)] 
=  0.
\end{align*}
Conversely, 
we see for
$(g,h) \in \Ker(\delta)$,
$$
\epsilon_{(1)}^{\lambda(\Lambda)}(g) + h = (\id - \lambda(\Lambda))(g^\prime)
\qquad
\text{ for some }
g^\prime \in  {\Bbb Z}_{\lambda(\Lambda)}.
$$
It follows  that
\begin{align*}
  & (\id - A_{\lambda(\Lambda)})(\lambda_\Lambda(1)^*(g^\prime) + g, g^\prime)\\
= & (\lambda_\Lambda(1)^*(g^\prime) + g - \lambda_\Lambda(1)^*(g^\prime), 
    g^\prime - \epsilon_{(1)}^{\lambda(\Lambda)}(\lambda_\Lambda(1)^*(g^\prime) + g) -  
       \sum_{\alpha=2}^{N} \epsilon_{(\alpha)}^{\lambda(\Lambda)} \circ\lambda_\Lambda(\alpha)^*(g^\prime))\\
= & (g, (\id - \lambda(\Lambda))(g^\prime) - \epsilon_{(1)}^{\lambda(\Lambda)}(g))
=  (g,h).
\end{align*}
Hence we have
$$
\Ker(\delta) = (\id - A_{\lambda(\Lambda)}) 
({\Bbb Z}_{\lambda({}_1\Lambda)} \oplus{\Bbb Z}_{\lambda(\Lambda)}).
$$

(ii)
Define a homomorphism
$$
\xi: h \in {\Bbb Z}_{\lambda(\Lambda)}
\longrightarrow
(\lambda_\Lambda(1)^*(h),h) \in {\Bbb Z}_{\lambda({}_1\Lambda)} 
\oplus{\Bbb Z}_{\lambda(\Lambda)}.
$$
An element
$ 
(g,h) \in {\Bbb Z}_{\lambda({}_1\Lambda)} \oplus{\Bbb Z}_{\lambda(\Lambda)}
$
belongs to
$
\Ker(\id - A_{\lambda(\Lambda)})
$
 if and only if 
 $
    (g,h) =   (\lambda_\Lambda(1)^*(h), 
   \epsilon_{(1)}^{\lambda(\Lambda)}(g) + \sum_{\alpha=2}^{N} 
   \epsilon_{(\alpha)}^{\lambda(\Lambda)} \circ\lambda_\Lambda(\alpha)^*(h)).
$
This condition is equivalent to the equalities:
$
g = \lambda_\Lambda(1)^*(h)
$ 
and
$ 
    h = \sum_{\alpha=1}^{N} \epsilon_{(\alpha)}^{\lambda(\Lambda)}
        \circ\lambda_\Lambda(\alpha)^*(h)
(= \lambda(\Lambda)(h)).
$
Since $\xi$ is clearly injective, 
it gives rise to an isomorphism:
$
\Ker (\id -\lambda(\Lambda))$ 
in 
$
{\Bbb Z}_{\lambda(\Lambda)}
\cong
\Ker (\id - A_{\lambda(\Lambda)})
$
 in 
$
{\Bbb Z}_{\lambda({}_1\Lambda}) \oplus{\Bbb Z}_{\lambda(\Lambda)}.
$

\end{pf}
Therefore we conclude:

\begin{thm}
For a $\lambda$-synchronizing subshift $\Lambda$,
the $\lambda$-synchronizing K-groups 
$K_0^\lambda (\Lambda), K_1^\lambda(\Lambda)$
 and
the $\lambda$-synchronizing Bowen-Franks groups 
$BF^0_\lambda(\Lambda),$
$ BF^1_\lambda(\Lambda)$
are all invariant under flow equivalence of $\lambda$-synchronizing subshifts.
\end{thm}
\begin{pf}
It is a direct consequence that the K-groups are invariant under flow
equivalence from  Theorem 4.4,  Proposition 5.11  Lemma 5.14 and Lemma 5.15.
As the Bowen-Franks groups are determined by the K-groups from the universal coefficient type theorem,
the groups $BF^*_\lambda(\Lambda)$ are also invariant under flow equivalence.
\end{pf}

\section{Examples}
{\bf 1. Sofic shifts}.

Let $\Lambda$ be an irreducible sofic shift 
and
 $\G_{F(\Lambda)}$
a finite directed labeled graph of the minimal left-resolving presentation of $\Lambda$,
that is called the left Fischer cover graph for $\Lambda$
(\cite{Fi}, cf. \cite{Kr84}, \cite{Kr87}, \cite{We}).
Let
${\frak L}_{\G_{F(\Lambda)}}$
be the $\lambda$-graph system
associated to the finite labeled graph ${\frak L}_{\G_{F(\Lambda)}}$
(see \cite[p. 20]{2002DocMath}).
Then the $\lambda$-synchronizing $\lambda$-graph system
$\LLL$ for the sofic shift $\Lambda$
is nothing but the $\lambda$-graph system
${\frak L}_{\G_{F(\Lambda)}}$.
Hence 
the $\lambda$-synchronizing $\lambda$-graph system
$\LLL$ for the sofic shift $\Lambda$
is regarded as the Fischer cover for the sofic shift $\Lambda$. 
Let $N$ be the number of the vertices of the graph $\G_{F(\Lambda)}$.
Let $\M_{F(\Lambda)}$
be the $N \times N$
symbolic matrix of the graph $\G_{F(\Lambda)}$.
Let  $A_{F(\Lambda)}$
be the $N \times N$ nonnegative matrix 
defined from $\M_{F(\Lambda)}$
by setting all symbols in each components  $\M_{F(\Lambda)}(i,j), i,j=1,\dots,N$
equal to $1$.
Then the $\lambda$-synchronizing $K$-groups and 
the $\lambda$-synchronizing Bowen-Franks groups 
are easily calculated as:
\begin{align*}
K_0^\lambda (\Lambda) & = {\Bbb Z}^N/(I_N - {}^t\negthinspace A_{F(\Lambda)}){\Bbb Z}^N,
\qquad 
 K_1^\lambda(\Lambda) = \Ker(I_N - {}^t\negthinspace A_{F(\Lambda)}) 
\quad \text{ in } {\Bbb Z}^N \\
\intertext{ and }
BF^0_\lambda(\Lambda)& = {\Bbb Z}^N/(I_N - A_{F(\Lambda)}){\Bbb Z}^N,
\qquad 
BF^1_\lambda(\Lambda) = \Ker(I_N - A_{F(\Lambda)}) \quad \text{ in } {\Bbb Z}^N.
\end{align*}
Therefore we have
\begin{prop}
Let $\Lambda$ be an irreducible sofic shift 
and
 $A_{F(\Lambda)}$
 the $N \times N$ edge matrix of its  left Fischer cover graph 
with entries in nonnegative integers.
 Then the abelian groups
 \begin{equation*}
 {\Bbb Z}^N/(I_N - A_{F(\Lambda)}){\Bbb Z}^N,\qquad
 \Ker(I_N - A_{F(\Lambda)}) \text{ in } {\Bbb Z}^N
\end{equation*}
are invariant under flow equivalence of subshifts (cf. \cite{FujO}).
\end{prop}
We note that 
$ {\Bbb Z}^N/(I_N - {}^t\negthinspace A_{F(\Lambda)}){\Bbb Z}^N$
is isomorphic to
${\Bbb Z}^N/(I_N - A_{F(\Lambda)}){\Bbb Z}^N$ 
and
$ \Ker(I_N - {}^t\negthinspace A_{F(\Lambda)}) \text{ in } {\Bbb Z}^N$
is isomorphic to
$
 \Ker(I_N - A_{F(\Lambda)}) \text{ in } {\Bbb Z}^N.
$

\medskip

{\bf 2. Dyck shifts.}

Let $N > 1$ be a fixed positive integer. 
Consider the Dyck shift $D_N$ 
with alphabet 
$\Sigma = \Sigma^- \sqcup \Sigma^+$
where
$\Sigma^- = \{ \alpha_1,\dots,\alpha_N \},
\Sigma^+ = \{ \beta_1,\dots,\beta_N \}.
$
The symbols 
$\alpha_i, \beta_i$
correspond to 
the brackets
$(_i,  )_i$
respectively.
The Dyck inverse monoid ${\Bbb D}_N$  has the relations
\begin{equation}
\alpha_i \beta_j
=
\begin{cases}
 {\bold 1} & \text{ if } i=j,\\
 0  & \text{ otherwise} 
\end{cases}  
\end{equation}
for 
$ i,j = 1,\dots,N$ (\cite{Kr}, \cite{Kr5}).
A word 
$\omega_1\cdots\omega_n$ 
of $\Sigma$
is admissible for $D_N$ 
precisely if
$
\prod_{m=1}^{n} \omega_m \ne 0.
$
For a word $\omega= \omega_1 \cdots \omega_n $ of $\Sigma,$ 
we denote by $\tilde{\omega}$ 
its reduced form, 
which is a word of 
$\Sigma \cup \{ 0, {\bold 1} \}$
obtained after the operation (6.1).
Hence a word $\omega$ of $\Sigma$
is forbidden for $D_N$ if 
and only if $\tilde{\omega} = 0$.

In \cite{KM2003}, 
a $\lambda$-graph system ${\frak L}^{Ch(D_N)}$ that presents $D_N$
has been introduced.
It is called the Cantor horizon $\lambda$-graph system for $D_N$
(cf. \cite{2007JOT}).
Let $\Sigma^N$ be the full $N$-shift 
$\{ 1,\dots,N \}^{\Bbb Z}$.
We denote by 
$B_l(D_N)$
and 
$B_l(\Sigma^N)$
 the set of admissible words of length 
$l$ of $D_N$
and that of 
$\Sigma^N$ respectively.
The vertices $V_l^{Ch(D_N)}$ 
of ${\frak L}^{Ch(D_N)}$
are given by the words of length $l$
consisting of the symbols of $\Sigma^+$.
That is, 
$$
V_l^{Ch(D_N)} = \{ \beta_{\mu_1}\cdots\beta_{\mu_l} \in B_l(D_N)
 \mid \mu_1\cdots\mu_l\in B_l(\Sigma^N) \}.
$$
The cardinal number of $V_l^{Ch(D_N)}$
 is $N^l$.
The mapping $\iota ( = \iota_{l,l+1}) :
V_{l+1}^{Ch(D_N)}\rightarrow V_l^{Ch(D_N)}$ 
deletes the rightmost of a word such as 
$$ 
\iota(\beta_{\mu_1}\cdots\beta_{\mu_{l+1}}) 
= \beta_{\mu_1}\cdots\beta_{\mu_l},
\qquad \beta_{\mu_1}\cdots\beta_{\mu_{l+1}} \in V_{l+1}^{Ch(D_N)}.
$$
There exists an edge labeled $\alpha_j$ from
the vertex 
$\beta_{\mu_1}\cdots\beta_{\mu_l}\in V_l^{Ch(D_N)}$ 
to the vertex
$\beta_{\mu_0}\beta_{\mu_1}\cdots\beta_{\mu_l}\in V_{l+1}^{Ch(D_N)}$
precisely if
$\mu_0 = j,$
and 
there exists an edge labeled $\beta_j$ from
$\beta_j\beta_{\mu_1}\cdots\beta_{\mu_{l-1}} \in V_l^{Ch(D_N)}$ 
to
$\beta_{\mu_1}\cdots\beta_{\mu_{l+1}} \in V_{l+1}^{Ch(D_N)}.$
The resulting labeled Bratteli diagram with  
$\iota$-map
is the Cantor horizon $\lambda$-graph system ${\frak L}^{Ch(D_N)}$ of $D_N$.
It is easy to see that each word of 
$V_l^{Ch(D_N)}$ is $l$-synchronizing in $D_N$
such that
$V_l^{Ch(D_N)}$ represent the all $l$-past equivalence classes of $D_N$.
Hence
we know that 
$V_l^{Ch(D_N)} = V^{\lambda(D_N)}_l.$
We then know that the canonical $\lambda$-synchronizing $\lambda$-graph system
${\frak L}^{\lambda(D_N)}$ is ${\frak L}^{Ch(D_N)}$.
\begin{prop}
The Dyck shift $D_N$ is $\lambda$-synchronizing, 
and
the $\lambda$-synchronizing $\lambda$-graph system 
${\frak L}^{\lambda(D_N)}$
is 
the Cantor horizon $\lambda$-graph system ${\frak L}^{Ch(D_N)}$.
\end{prop}
The K-groups of the $\lambda$-graph system ${\frak L}^{Ch(D_N)}$
have been  computed  in \cite{KM2003} and \cite{2007JOT}
so that we have  
$$
K^\lambda_0(D_N) 
\cong {\Bbb Z}/N{\Bbb Z} \oplus C({\frak K},{\Bbb Z}), \qquad 
K^\lambda_1(D_N) \cong 0,
$$
where 
$
C({\frak K},{\Bbb Z})
$
denotes the abelian group of all ${\Bbb Z}$-valued continuous functions 
on a Cantor discontinuum ${\frak K}$.
By Theorem 5.16, we have
\begin{prop}
For 
the Dyck shifts $D_N, D_{N'}$
with $N, N'\ge 2$,
$D_N$ is flow equivalent to 
$D_{N'}$
if and only if $N = N'$.
\end{prop}

\medskip

{\bf 3. Topological Markov Dyck shifts.}

We will state a generalization of 
the  Dyck shifts.
Let $A =[A(i,j)]_{i,j=1,\dots,N}$
be an $N\times N$ matrix with entries in $\{0,1\}$.
Consider the Dyck inverse monoid 
for the alphabet
$\Sigma = \Sigma^- \sqcup \Sigma^+$
where
$\Sigma^- = \{ \alpha_1,\cdots,\alpha_N \},
\Sigma^+ = \{ \beta_1,\cdots,\beta_N \}
$
as in the above.
Let ${\cal O}_A$
be  the Cuntz-Krieger algebra  for the matrix $A$
that is the universal $C^*$-algebra generated by 
$N$ partial isometries $t_1,\dots,t_N$ 
subject to the following relations:
$$
\sum_{j=1}^N t_j t_j^* = 1, 
\qquad
t_i^* t_i = \sum_{j=1}^N A(i,j) t_jt_j^* \quad \text{ for } i= 1,\dots,N
$$
(\cite{CK}).
Define a correspondence 
$\varphi_A :\Sigma \longrightarrow \{t_i^*, t_i \mid i=1,\dots,N\}$
by setting
$$ 
\varphi_A(\alpha_i) = t_i^*,\qquad 
\varphi_A(\beta_i) = t_i,  \qquad i=1,\dots,N.
$$
We denote by $\Sigma^*$ the set of all words 
$\gamma_1\cdots \gamma_n$ of elements of $\Sigma$.
Define the set
$$
{\frak F}_A = \{ \gamma_1\cdots \gamma_n \in \Sigma^* \mid
\varphi_A(\gamma_1)\cdots \varphi_A( \gamma_n) = 0 \text{ in } {\cal O}_A \}.
$$
Let $D_A$ be the subshift over $\Sigma$ whose forbidden words are 
${\frak F}_A.$
The subshift is called the topological Markov Dyck shift defined by $A$.
These kinds of  subshifts have first appeared  
in semigroup setting (\cite{HIK})
and 
in more general setting (\cite{KM2003}) without using $C^*$-algebras.
If all entries of $A$  are $1$, 
the subshift is nothing but the Dyck shift $D_N$
with $2N$ bracket, because  
the partial isometries 
$\{ \varphi_A(\alpha_i), \varphi(\beta_i) \mid i=1,\dots,N\}$
yield the  Dyck inverse monoid.
We note the fact that 
 $\alpha_i \beta_j\in {\frak F}_A$  if $i\ne j$,
 and
 $\alpha_{i_n}\cdots \alpha_{i_1} \in {\frak F}_A$  
if and only if 
$\beta_{i_1}\cdots \beta_{i_n} \in {\frak F}_A$.
Consider the following  two subsystems of $D_A$
\begin{align*}
D_A^+ & = \{ {(\gamma_i)}_{i \in \Bbb Z} \in D_A \mid
\gamma_i \in \Sigma^+, i \in \Bbb Z \},\\ 
D_A^- & = \{ {(\gamma_i)}_{i \in \Bbb Z} \in D_A \mid
\gamma_i \in \Sigma^-, i \in \Bbb Z \}. 
\end{align*}
The subshift  
$D_A^+$ is identified with the topological Markov shift 
$$
\Lambda_A = \{ {(x_i)}_{i \in \Bbb Z}\in \{ 1,\dots,N \}^{\Bbb Z}
 \mid A(x_i,x_{i+1}) = 1, i \in \Bbb Z \}
$$ 
defined by the matrix $A$ and similarly
$D_A^-$ is identified with 
the topological Markov shift $\Lambda_{{}^t\negthinspace A}$ 
defined by the transposed matrix ${}^t\negthinspace A$ of $A$.
Hence  $D_A$ contains the both topological Markov shifts 
$\Lambda_A$ and $\Lambda_{{}^t\negthinspace A}$ that do not intersect each other.
If $A$ satisfies condition (I) in the sense of Cuntz-Krieger \cite{CK},
the subshift $D_A$ is not sofic (\cite[Proposition 2.1]{2010MS}).
We may define a $\lambda$-graph system 
${\frak L}^{Ch(D_A)}$
called  the Cantor horizon $\lambda$-graph system for $D_A$ 
(\cite{2010MS}).
We denote by 
$B_l(D_A)$
and 
$B_l(\Lambda_A)$
 the set of admissible words of length 
$l$ of $D_A$
and that of 
$\Lambda_A$ respectively.
The vertices $V_l^{Ch(D_A)}$ of ${\frak L}^{Ch(D_A)}$ at level $l$
are given by the admissible words of length $l$
consisting of the symbols of $\Sigma^+$.
They are $l$-synchronizing words of $D_A$
such that
the $l$-past equivalence classes of them
coincide with the $l$-past equivalence classes of 
the set of  $l$-synchronizing words of $D_A$.
Hence $V_l^{Ch(D_A)} = V^{\lambda(D_A)}_l$.
Since $V_l^{Ch(D_A)}$ 
is identified with $B_l(\Lambda_A)$,
we may write $V_l^{Ch(D_A)}$ as
$$
V_l^{Ch(D_A)} = \{ v^l_{\mu_1 \cdots \mu_l}
 \mid \mu_1\cdots\mu_l\in B_l(\Lambda_A) \}.
$$
The mapping $\iota ( = \iota_{l,l+1}) :V_{l+1}^{Ch(D_A)} \rightarrow V_l^{Ch(D_A)}$ 
is defined by deleting the rightmost symbol of a corresponding word such as 
$$ 
\iota( v^{l+1}_{\mu_1\cdots \mu_{l+1}}) 
= v^l_{\mu_1 \cdots \mu_l},
\qquad v^{l+1}_{\mu_1 \cdots \mu_{l+1}} \in V_{l+1}^{Ch(D_A)}.
$$
There exists an edge labeled $\alpha_j$ from
$v^l_{\mu_1 \cdots \mu_l}\in V_l^{Ch(D_A)}$ 
to
$v^{l+1}_{\mu_0\mu_1 \cdots \mu_l}\in V_{l+1}^{Ch(D_A)}$
precisely if
$\mu_0 = j,$
and 
there exists an edge labeled $\beta_j$ from
$v^l_{j\mu_1 \cdots \mu_{l-1}}
\in V_l^{Ch(D_A)}$ 
to
$v^{l+1}_{\mu_1 \cdots \mu_{l+1}}\in V_{l+1}^{Ch(D_A)}.$
The resulting labeled Bratteli diagram with 
$\iota$-map
becomes a $\lambda$-graph system
written 
${\frak L}^{Ch(D_A)}$
that presents $D_A$.
It is called 
the Cantor horizon $\lambda$-graph system 
for $D_A$.

\begin{prop}
The subshift $D_A$ is $\lambda$-synchronizing,
and 
 the $\lambda$-synchronizing $\lambda$-graph system
${\frak L}^{\lambda(D_A)}$ is
the Cantor horizon $\lambda$-graph system
${\frak L}^{Ch(D_A)}$.
\end{prop}
By \cite[Lemma 2.5]{2010MS}
 if $A$ satisfies condition (I) in the sense of Cuntz-Krieger \cite{CK},
the $\lambda$-graph system ${\frak L}^{Ch(D_A)}$ 
satisfies $\lambda$-condition (I).
 If $A$ is irreducible, it is $\lambda$-irreducible.
In this case we have that 
 the $C^*$-algebra 
${\cal O}_{{\frak L}^{\lambda(D_A)}}$
associated with ${\frak L}^{\lambda(D_A)}$
is simple, purely infinite
(cf.\cite{2010MS}).

\medskip

One knows that 
$\beta$-shifts for $1 <\beta\in {\Bbb R}$,
a synchronizing counter shift named as the context free shift in \cite[Example 1.2.9]{LM},
Motzkin shifts and
the Morse shift
are all $\lambda$-synchronizing.
Their proofs are essentially seen in the papers
 \cite{KMW}, \cite{KM2010}, \cite{1999JOT}, \cite{2004MZ}
respectively.

\medskip

We  study $C^*$-algebras associated with $\lambda$-synchronizing $\lambda$-graph systems
in \cite{Ma2011}.


\end{document}